\documentclass{article}

\title{Multiplicity of solutions of a zero mass nonlinear equation on 
  a Riemannian manifold}

\author{D.~Visetti\thanks{Department of Mathematics, University 
    of Trento, Italy; e-mail: \texttt{visetti@science.unitn.it}.}}

\RequirePackage[italian,english]{babel} 
\RequirePackage[dvips]{graphicx}
\RequirePackage{amsmath}
\RequirePackage{amssymb}
\RequirePackage{amsthm}

\numberwithin{equation}{section}

\begin{document}
\selectlanguage{english}

\newcommand{\NN}{{\mathbb N}}
\newcommand{\RR}{{\mathbb R}}
\newcommand{\ZZ}{{\mathbb Z}}
\newcommand{\QQ}{{\mathbb Q}}
\newcommand{\CC}{{\mathbb C}}

\newtheorem{theorem}{Theorem}[section]
\newtheorem{proposition}[theorem]{Proposition}
\newtheorem{lemma}[theorem]{Lemma}
\newtheorem{corollary}[theorem]{Corollary}
\newtheorem{remark}[theorem]{Remark}
\newtheorem{definition}[theorem]{Definition}

\maketitle


\section{Introduction}

In this paper we are interested in the relation between the number 
of solutions of a nonlinear equation on a Riemannian manifold and 
the topology of the manifold itself.

Let $(M,g)$ be a compact, connected, orientable, boundaryless 
Riemannian manifold of class $C^\infty$ with Riemannian metric $g$.  
Let $\dim(M)=n\geq 3$.

We consider the problem
\begin{equation}
-\epsilon^2\Delta u=f'(u)
\label{eq}
\end{equation}
with $u\in H_1^2(M)$.

As it has been pointed out in \cite{BL} problem (\ref{eq}) admits 
solutions on $\RR^n$ if $f'(0)<0$, while there are no solutions 
if $f'(0)>0$.  The limiting case $f'(0)=0$, i.e. the ``zero mass'' case, 
depends on the structure of $f$.  Berestycki and Lions proved the 
existence of ground state solutions if $f(u)$ behaves as $|u|^p$ 
for $u$ large and $|u|^q$ for $u$ small, with $p$ and $q$ respectively 
super and sub-critical.   In \cite{BLzero} they proved also the 
existence of infinitely many bound state solutions.

Problem (\ref{eq}) has been studied also in \cite{BM}, where existence 
and non existence results have been given on an exterior domain in 
$\RR^n$.

The problem of the multiplicity of solutions of a nonlinear elliptic 
equation on a Riemannian manifold has been studied 
in \cite{bbm}, where the authors consider an equation with sub-critical 
growth.

The effect of the domain shape on the number of positive solutions 
of some semilinear elliptic problems has been widely studied.  Here we 
only mention \cite{bahri-coron}, \cite{D}, \cite{bc}, \cite{bcp} and 
\cite{BC}.

Let $f:\RR\to\RR$ be an even function such that:
\begin{itemize}
\item[(f1)] $0<\mu f(s)\leq f'(s)s<f''(s)s^2$ for any $s\neq 0$ and for 
  some $\mu>2$;
\item[(f2)] $f(0)=f'(0)=f''(0)=0$ and there exist positive constants 
  $c_0,c_1,p,q$ with $2<p<2^*<q$ such that
  \begin{eqnarray}
  f(s) & \geq & \left\{ \begin{array}{ll}
                c_0|s|^p & \mbox{for }|s|\geq 1 \\
                c_0|s|^q & \mbox{for }|s|\leq 1
                \end{array}\right. \label{f}            \\
  f''(s) & \geq & \left\{ \begin{array}{ll}
                c_1|s|^{p-2} & \mbox{for }|s|\geq 1 \\
                c_1|s|^{q-2} & \mbox{for }|s|\leq 1 
                \end{array}\right. \label{f''}
  \end{eqnarray}
\end{itemize}
We denote by $\mathrm{cat}\,(M)$ the Ljusternik-Schnirelmann category 
of $M$ and by $\mathcal{P}_t(M)$ the Poincar\'e polynomial of $M$.

Our main results are the following:

\begin{theorem}
\label{trm-gen1}
  For $\epsilon>0$ sufficiently small, equation (\ref{eq}) has at least 
  $\mathrm{cat}\,(M)+1$ solutions in $H_1^2(M)$.
\end{theorem}

\begin{theorem}
\label{trm-morse}
  If for $\epsilon>0$ sufficiently small the solutions of equation 
  (\ref{eq}) are non-degenerate, then there are at least 
  $2\mathcal{P}_1(M)-1$ solutions.
\end{theorem}

\section{Notation and preliminary results}
\label{sec:notation}

We denote by $B(0,R)$ the ball in $\RR^n$ of centre $0$ and radius $R$ 
and by $B_g(x,R)$ the ball in $M$ of centre $x$ and radius $R$.

We define a smooth real function $\chi_R$ on $\RR^+$ such that
\begin{equation}
\chi_R(t) = \left\{ \begin{array}{ll}
  1 & \mbox{if }0\leq t\leq\frac{R}{2} \\
  0 & \mbox{if }t\geq R
  \end{array} \right.
\label{cut-off}
\end{equation}
and $|\chi_R'(t)|\leq\frac{\chi_0}{R}$, with $\chi_0$ positive constant.


We recall some definitions and results about compact connected 
Riemannian manifolds of class $C^\infty$ (see for example \cite{H}).

\begin{remark}
\label{rmk-exp}
  On the tangent bundle $TM$ of $M$ the exponential map $\exp:TM\to M$ 
  is defined.  This map has the following properties:
  \begin{itemize}
  \item[(i)] $\exp$ is of class $C^\infty$;
  \item[(ii)] there exists a constant $R>0$ such that
    $$
    \exp_x|_{B(0,R)}:B(0,R)\to B_g(x,R)
    $$
    is a diffeomorphism for all $x\in M$.
  \end{itemize}
\end{remark}

It is possible to choose an atlas $\mathcal{C}$ on $M$, whose charts are 
given by the exponential map (normal coordinates).  We denote by 
$\{\psi_{C}\}_{C\in\mathcal{C}}$ a partition of unity subordinate to 
the atlas $\mathcal{C}$.   Let $g_{x_0}$ be the Riemannian metric 
in the normal coordinates of the map $\exp_{x_0}$.

For any $u\in H_1^2(M)$ we have that:
$$
\begin{aligned}
  &\int_M |\nabla u(x)|^2_g d\mu_g
   = \sum_{C\in\mathcal{C}}\int_C\psi_C(x)|\nabla u(x)|^2_g d\mu_g \\
  &= \sum_{C\in\mathcal{C}}\int_{B(0,R)}\psi_C(\exp_{x_C}(z))\,
   g^{ij}_{x_C}(z)\frac{\partial u(\exp_{x_C}(z))}{\partial z_i}
   \frac{\partial u(\exp_{x_C}(z))}{\partial z_j}\,|g_{x_C}(z)|^
   \frac{1}{2}\, dz\, ,
\end{aligned}
$$
where Einstein notation is adopted, that is
$$
g^{ij}z_iz_j=\sum_{i,j=1}^ng^{ij}z_iz_j\, ,
$$
$(g^{ij}_{x_0}(z))$ is the inverse matrix of $g_{x_0}(z)$ and 
$|g_{x_0}(z)|=\det(g_{x_0}(z))$.   In particular we have that 
$g_{x_0}(0)=\mathrm{Id}$.  
A similar relation holds for the integration of $|u(x)|^p$.  For 
convenience we will also write for all $x_0\in M$ and $z,\xi\in T_{x_0}M$
\begin{equation}
\label{normag}
|\xi|^2_{g_{x_0}(z)} = g^{ij}_{x_0}(z)\xi_i\xi_j\, .
\end{equation}

\begin{remark}
\label{rmk-metric}
  Since $M$ is compact, there are two strictly positive constants $h$ 
  and $H$ such that for all $x\in M$ and all $z\in T_xM$
  $$
  h|z|^2\leq g_x(z,z)\leq H|z|^2\, ,
  $$
  where $|\cdot|$ is the standard metric in $\RR^n$.  Hence there holds
  $$
  h^n\leq |g_x(z)|\leq H^n\, .
  $$
\end{remark}

We are going to find the solutions of (\ref{eq}) as critical points of the 
functional $J_\epsilon:H_1^2(M)\to\RR$, defined by
\begin{equation}
J_\epsilon(u) = \frac{\epsilon^2}{2\epsilon^n} \int_M |\nabla u(x)|^2_g
  \, d\mu_g - \frac{1}{\epsilon^n} \int_M f(u(x))\, d\mu_g\, ,
\label{Jepsilon}
\end{equation}
constrained on the Nehari manifold
\begin{equation}
\mathcal{N}_\epsilon = \left\{ u\in H_1^2(M) \;\big\vert\;u\neq 0
  \mbox{ and }\int_M\epsilon^2|\nabla u|^2_g\, d\mu_g=\int_Mf'(u)u\, 
  d\mu_g\right\}.
\label{Nepsilon}
\end{equation}

Let $\mathcal{D}^{1,2}(\RR^n)$ be the completion of $C_0^\infty(\RR^n)$ 
with respect to the norm
$$
\|v\|^2_{\mathcal{D}^{1,2}(\RR^n)}=\int_{\RR^n}
|\nabla v(z)|^2dz\, .
$$
We consider also the following functional $J:\mathcal{D}^{1,2}(\RR^n)
\to\RR$ defined by
\begin{equation}
J(v) := \int_{\RR^n}\left(\frac{1}{2}|\nabla v(x)|^2-f(v(x))\right)dx
\label{J}
\end{equation}
and the associated Nehari manifold
\begin{equation}
\mathcal{N} = \left\{ v\in\mathcal{D}^{1,2}(\RR^n)\;\big\vert\;
  v\neq 0\mbox{ and }\int_{\RR^n}|\nabla v(x)|^2\, dx=\int_{\RR^n}
  f'(u)u\, dx\right\}\, .
\label{N}
\end{equation}

The functionals $J_\epsilon$ and $J$ are $C^2$ respectively on 
$H_1^2(M)$ and on $\mathcal{D}^{1,2}(\RR^n)$.  In fact, we have:

\begin{lemma}
\label{lmm-FC2}
  \begin{itemize}
  \item[(i)] The functional $F_{\epsilon,M}:L^p(M)\to\RR$, 
    defined by
    \begin{equation}
    \label{FepsilonM}
    F_{\epsilon,M}(u) := \frac{1}{\epsilon^n}\int_M f(u(x))\, d\mu_g
    \end{equation}
    is of class $C^2$ and
    \begin{eqnarray*}
    F_{\epsilon,M}'(u_0)u_1     & = & \frac{1}{\epsilon^n}\int_M 
                                      f'(u_0(x))u_1(x)\, d\mu_g  \\
    F_{\epsilon,M}''(u_0)u_1u_2 & = & \frac{1}{\epsilon^n}\int_M 
                                      f''(u_0(x))u_1(x)u_2(x)\, d\mu_g
    \end{eqnarray*}
  \item[(ii)] The functional $F:L^{2^*}(\RR^n)\to\RR$ defined by
    \begin{equation}
    \label{F}
    F(v) := \int_{\RR^n} f(v(z))\, dz
    \end{equation}
    is of class $C^2$ and
    \begin{eqnarray*}
    F'(v_0)v_1     & = & \int_{\RR^n} f'(v_0(z))v_1(z)\, dz \\
    F''(v_0)v_1v_2 & = & \int_{\RR^n} f''(v_0(z))v_1(z)v_2(z)\, dz
    \end{eqnarray*}
  \end{itemize}
\end{lemma}

The proof of this lemma is analogous to the proof of Lemma 2.7 in 
\cite{BM}.

We also have the following lemma:

\begin{lemma}
\label{lmm-Ftilde}
  The functionals $\widetilde F_{\epsilon,M}:L^p(M)\to\RR$, defined by
  \begin{equation}
  \label{FtildeepsilonM}
  \widetilde F_{\epsilon,M}(u) := \frac{1}{\epsilon^n}\int_M \left[
    \frac{1}{2}f'(u(x))u(x)-f(u(x))\right] d\mu_g
  \end{equation}
  and $\widetilde F_\Omega:L^{2^*}(\Omega)\to\RR$ defined by
  \begin{equation}
  \label{FtildeOmega}
  \widetilde F_\Omega(v) := \int_{\Omega} \left[\frac{1}{2}f'(v(z))v(z)-
    f(v(z))\right] dz
  \end{equation}
  are strongly continuous.
\end{lemma}

We write
\begin{equation}
m(J) := \inf\{J(v)\;\vert\;v\in\mathcal{N}\}\, .
\label{m(J)}
\end{equation}

There exists a positive spherically symmetric and decreasing with $|z|$ 
solution $U\in\mathcal{D}^{1,2}(\RR^n)$ of
\begin{equation}
-\Delta U=f'(U)\quad\mbox{in}\quad\RR^n\, ,
\label{eqRn}
\end{equation}
such that $J(U)=m(J)$ (see \cite{BL} and \cite{BM}).

The function $U_\epsilon(z)=U\left(\frac{z}{\epsilon}\right)$ is 
solution of
$$
-\epsilon^2\Delta U_\epsilon=f(U_\epsilon)\, .
$$

For any $\delta>0$ we consider the subset of $\mathcal{N}_\epsilon$
\begin{equation}
\Sigma_{\epsilon,\delta} := \{u\in\mathcal{N}_\epsilon\;\vert\;
  J_\epsilon(u)<m(J)+\delta\}\, .
\label{Sigmaepsilondelta}
\end{equation}

We recall now the definition of Palais-Smale condition:

\begin{definition}
  Let $J$ be a $C^1$ functional on a Banach space $X$.  A sequence 
  $\{u_m\}$ in $X$ is a Palais-Smale sequence for $J$ if $|J(u_m)|
  \leq c$, uniformly in $m$, while $J'(u_m)\to 0$ strongly, as $m\to
  \infty$.   We say that $J$ satisfies the Palais-Smale condition 
  ((PS) condition) if any Palais-Smale sequence has a convergent 
  subsequence.
\end{definition}

\section{Ideas of the proof for the category theory result}
\label{sec:ideas}

We recall the definition of Ljusternik-Schnirelmann category (see 
\cite{LS}).

\begin{definition}
  Let $M$ be a topological space and consider a closed subset $A\subset 
  M$.   We say that $A$ has category $k$ relative to $M$ 
  ($\mathrm{cat}_M(A)=k$) if $A$ is covered by $k$ closed sets $A_j$, 
  $1\leq j\leq k$, which are contractible in $M$ and if $k$ is minimal 
  with this property.   If no such finite covering exists, we let 
  $\mathrm{cat}_M(A)=\infty$.  If $A=M$, we write $\mathrm{cat}_M(M)=
  \mathrm{cat}(M)$.
\end{definition}

\begin{remark}
\label{rmk-cat}
  Let $M_1$ and $M_2$ be topological spaces.   If $g_1:M_1\to M_2$ and 
  $g_2:M_2\to M_1$ are continuous operators such that $g_2\circ g_1$ is 
  homotopic to the identity on $M_1$, then $\mathrm{cat}\,(M_1)\leq
  \mathrm{cat}\,(M_2)$ (see \cite{bc}).
\end{remark}

Using the notation in the previous section, Theorem~\ref{trm-gen1} can 
be stated more precisely like this:

\begin{theorem}
\label{trm1}
  There exists $\delta_0\in(0,m(J))$ such that for any $\delta\in(0,
  \delta_0)$ there exists $\epsilon_0=\epsilon_0(\delta)>0$ and for any 
  $\epsilon\in(0,\epsilon_0)$ the functional $J_\epsilon$ has at least 
  $\mathrm{cat}\,(M)$ critical points $u\in H^1_2(M)$ satisfying 
  $J_\epsilon(u)<m(J)+\delta$ and at least one critical point 
  with $J_\epsilon(u)>m(J)+\delta$.
\end{theorem}

This theorem is a consequence of the following classical result (see 
for example \cite{bcp}):

\begin{theorem}
\label{trm-cat}
  Let $J$ be a $C^1$ real functional on a complete $C^{1,1}$ 
  submanifold $N$ of a Banach space. If $J$ is bounded below and 
  satisfies the (PS) condition then it has at least $\mathrm{cat}(J^d)$ 
  critical points in $J^d$, where $J^d:=\{u\in N : J(u)<d\}$, 
  and at least one critical point $u \not\in J^d$.
\end{theorem}

More precisely, Theorem~\ref{trm1} follows from the previous theorem, 
Remark~\ref{rmk-cat} and the following proposition:

\begin{proposition}
\label{prp-cat}
  There exists $\delta_0\in(0,m(J))$ such that for any $\delta\in(0,
  \delta_0)$ there exists $\epsilon_0=\epsilon_0(\delta)>0$ and for any 
  $\epsilon\in(0,\epsilon_0)$ we have
  $$
  \mathrm{cat}\,(M)\leq\mathrm{cat}\,(\Sigma_{\epsilon,\delta})\, .
  $$
\end{proposition}

In order to prove this we will present two suitable functions $g_1$ 
and $g_2$.

By the embedding theorem, we assume that $M$ is embedded in $\RR^N$, 
with $N\geq 2n$.

\begin{definition}
\label{def-r(M)}
  We define the radius of topological invariance $r(M)$ of $M$ as
  $$
  r(M) := \sup\{\rho>0\ |\ \mathrm{cat}\,(M_\rho)=\mathrm{cat}\,(M)\}\, ,
  $$
  where $M_\rho:=\{z\in\RR^N\ |\ d(z,M)<\rho\}$.
\end{definition}

We can now show a function $\phi_\epsilon:M\to\Sigma_{\epsilon,\delta}$ 
and a function $\beta:\Sigma_{\epsilon,\delta}\to M_r$, with $0<r<r(M)$ 
such that
\begin{equation}
I_\epsilon := \beta\circ\phi_\epsilon:M\to M_r
\label{Iepsilon}
\end{equation}
is well defined and homotopic to the identity on $M$.

\section{The function $\phi_\epsilon$}

Next lemma presents some properties of the Nehari manifold.

\begin{lemma}
\label{lmm-tepsilonu}
  \begin{itemize}
  \item[(i)] The set $\mathcal{N}_\epsilon$ (resp. $\mathcal{N}$) 
    is a $C^1$ manifold.
  \item[(ii)] For all not constant $u\in H_1^2(M)$ (for all $v\in
    \mathcal{D}^{1,2}(\RR^n)$, $v\not\equiv0$), there exists a unique 
    $t_\epsilon(u)>0$ ($t(v)>0$) such that $t_\epsilon(u)u\in
    \mathcal{N}_\epsilon$ ($t(v)v\in\mathcal{N}$) and $J_\epsilon
    (t_\epsilon(u)u)$ ($J(t(v)v)$) is the maximum value of 
    $J_\epsilon(tu)$ ($J(tv)$) for $t\geq 0$.
  \item[(iii)] The dependence of $t_\epsilon(u)$ on $u$ (of $t(v)$ on 
    $v$) is $C^1$.
  \end{itemize}
\end{lemma}

For the proof see Lemma 3.1 in \cite{BM}.

Let $U$ be the function defined in Section~\ref{sec:notation}.  We write
$$
\widetilde U_\frac{R}{\epsilon} = U(z) \quad \mbox{with} \quad z\in\RR^n 
  \quad \mbox{such that} \quad |z|=\frac{R}{\epsilon}\, .
$$
For any $x_0\in M$ and $\epsilon>0$, we consider the function on $M$
\begin{equation}
W_{x_0,\epsilon}(x) := \left\{ \begin{array}{ll}
    U_\epsilon(\exp_{x_0}^{-1}(x))-\widetilde U_\frac{R}{\epsilon} & 
    \mbox{if }x\in B_g(x_0,R)\, ,                                  \\
    0 & \mbox{otherwise,}
  \end{array} \right.
\label{Wx0epsilon}
\end{equation}
where $R$ is chosen as in Remark~\ref{rmk-exp} $(ii)$.

The function $W_{x_0,\epsilon}$ is in $H_1^2(M)$ and is not identically 
zero.  Then, by the previous lemma, we can define
\begin{equation}
\begin{array}{ccccc}
  \phi_\epsilon & : & M & \longrightarrow & \mathcal{N}_\epsilon        \\
  &  & x_0 & \longmapsto & t_\epsilon(W_{x_0,\epsilon}(x))W_{x_0,
    \epsilon}(x)\, .
\end{array}
\label{phiepsilon}
\end{equation}

The choice of the function $\phi_\epsilon$ different from the one in 
\cite{bbm} has been made for the function $U$ can be not in $L^2(\RR^n)$.

\begin{proposition}
  For any $\epsilon>0$ the map $\phi_\epsilon:M\to
  \mathcal{N}_\epsilon$ is continuous.  For any $\delta>0$ there 
  exists $\epsilon_0>0$ such that if $\epsilon<\epsilon_0$
  $$
  \phi_\epsilon(x_0)\in\Sigma_{\epsilon,\delta}
  $$
  for all $x_0\in M$.
\end{proposition}

\begin{proof}
  \textbf{(I)} \emph{The map $\phi_\epsilon:M\to\mathcal{N}_\epsilon$ 
  is continuous.}

  By Lemma~\ref{lmm-tepsilonu} $(iii)$, it is enough to prove that
  $$
  \lim_{k\to\infty}\|W_{x_k,\epsilon}-W_{\hat x,\epsilon}\|_{H^1_2(M)}=0\, .
  $$
  for any sequence $\{x_k\}$ in $M$, converging to $\hat x$.
  
  We choose a finite atlas $\mathcal{C}$ for $M$, which contains the 
  chart $C=B_g(\hat x,R)$.   The functions $W_{x_k,\epsilon}$ and 
  $W_{\hat x,\epsilon}$ have support respectively on $B_g(x_k,R)$ 
  and on $B_g(\hat x,R)$.   Since $x_k\to\hat x$ the set $Z_k=[B_g(x_k,
  R)\setminus B_g(\hat x,R)]\cup[B_g(\hat x,R)\setminus B_g(x_k,R)]$ 
  is such that $\mu_g(Z_k)\to 0$ as $k\to\infty$. Then we have
  $$
  \int_{Z_k} \left| \nabla \left( W_{x_k,\epsilon}(x)-W_{\hat x,
    \epsilon}(x) \right)\right|^2_g d\mu_g \to 0 \qquad \mbox{as }
    k\to\infty\, .
  $$
  We still have to check the integral on $B_g(x_k,R)\cap B_g(\hat x,R)$.  
  We write $A_k=\exp_{\hat x}^{-1}(B_g(x_k,R)\cap B_g(\hat x,R))$ and 
  $\eta_k(z)=\exp_{x_k}^{-1}(\exp_{\hat x}(z))$
  $$
  \begin{aligned}
  \int_{\exp_{\hat x}(A_k)} \hspace{-1cm}\left| \nabla \left[ 
      W_{x_k,\epsilon}(x)-W_{\hat x,\epsilon}(x) \right]
      \right|^2_g d\mu_g
    &= \int_{A_k} \hspace{-0.2cm}\left| \nabla \left[ 
      U_\epsilon(\eta_k(z))-U_\epsilon(z) \right] \right|^2_
      {g_{\hat x}(z)} |g_{\hat x}(z)|^\frac{1}{2} dz           \\
    &\leq \frac{H^\frac{n}{2}}{h} \int_{A_k} \left| \nabla 
      \left[ U_\epsilon(\eta_k(z))-U_\epsilon(z) \right] 
      \right|^2 dz\, .
  \end{aligned}
  $$
  Since $\eta_k(z)$ tends point-wise to $z$ and $\nabla U_\epsilon$ 
  is continuous, $\left|\nabla[U_\epsilon(\eta_k(z))-U_\epsilon(z)]
  \right|^2$ tends pointwise to zero.  Applying Lebesgue theorem, 
  we obtain that
  $$
  \int_{M} \left| \nabla \left[ W_{x_k,\epsilon}(x)-W_{\hat x,\epsilon}
    (x) \right]\right|^2_gd\mu_g \to 0\, .
  $$
  In an analogous way we have that $\| W_{x_k,\epsilon}-W_{\hat x,
  \epsilon}\|^2_{L^2(M)}$ tends to zero.    

  \noindent\textbf{(II)} \emph{The limit of $\frac{\epsilon^2}{\epsilon^n}
  \int_M \left| \nabla W_{x_0,\epsilon}(x) \right|^2_gd\mu_g$ is 
  $\| U\|^2_{\mathcal{D}^{1,2}(\RR^n)}$}.

  To prove the second statementof this proposition, first we show that
  \begin{equation}
  \lim_{\epsilon\to 0} \frac{\epsilon^2}{\epsilon^n} \int_M \left| 
    \nabla W_{x_0,\epsilon}(x) \right|^2_gd\mu_g = \| U\|^2_{\mathcal
    {D}^{1,2}(\RR^n)}
  \label{limW1}
  \end{equation}
  uniformly with respect to $x_0\in M$.

  We evaluate the following:
  \begin{eqnarray*}
  \lefteqn{\left|\frac{\epsilon^2}{\epsilon^n} \int_M \left| 
      \nabla W_{x_0,\epsilon} \right|^2_g d\mu_g - \int_{\RR^n} 
      |\nabla U|^2\, dz \right|}                                  \\
    & = & \left|\frac{\epsilon^2}{\epsilon^n} \int_{B_g(x_0,R)} 
      \left|\nabla\left[U_\epsilon(\exp_{x_0}^{-1}(x))\right]
      \right|^2_g d\mu_g - \int_{\RR^n} |\nabla U|^2\, dz \right| \\
    & = & \left|\frac{\epsilon^2}{\epsilon^n} \int_{B(0,R)} 
      \left|\nabla U_\epsilon(z) \right|^2_{g_{x_0}(z)} |g_{x_0}
      (z)|^\frac{1}{2} dz - \int_{\RR^n} |\nabla U|^2\, dz 
      \right|\, .
  \end{eqnarray*}
  Changing variables, we obtain
  $$
  \left| \int_{\RR^n} \left(\chi_{B\left(0,\frac{R}{\epsilon}\right)}(z)
    g^{ij}_{x_0}(\epsilon z)|g_{x_0}(\epsilon z)|^\frac{1}{2}-
    \delta^{ij} \right) \frac{\partial U}{\partial z_i}\frac{\partial U}
    {\partial z_j} dz \right|,
  $$
  where $\chi_{B\left(0,\frac{R}{\epsilon}\right)}(z)$ denotes the 
  characteristic function of the set $B\left(0,\frac{R}{\epsilon}\right)$ 
  and where $\delta^{ij}$ is the Kronecker delta 
  (it takes value $0$ for $i\neq j$ and $1$ for $i=j$).  The previous 
  integral is bounded from above by the following sum
  $$
  \begin{aligned}
  \lefteqn{\left| \int_{B(0,T)} \left(\chi_{B\left(0,\frac{R}{\epsilon}
      \right)}(z)g^{ij}_{x_0}(\epsilon z)|g_{x_0}(\epsilon z)|^
      \frac{1}{2}-\delta^{ij} \right) \frac{\partial U}{\partial 
      z_i}\frac{\partial U}{\partial z_j} dz \right|}                  \\
    &+ \left| \int_{\RR^n\setminus B(0,T)} \left( \chi_{B\left(0,
      \frac{R}{\epsilon}\right)}(z)g^{ij}_{x_0}(\epsilon z)
      |g_{x_0}(\epsilon z)|^\frac{1}{2}-\delta^{ij} \right) 
      \frac{\partial U}{\partial z_i}\frac{\partial U}{\partial z_j} 
      dz \right|,
  \end{aligned}
  $$
  with $T>0$.  It is easy to see that the second addendum vanishes as 
  $T\to\infty$.  As regards the first addendum, fixed $T$, by compactness of 
  the manifold $M$ and regularity of the Riemannian metric $g$ the limit
  $$
  \lim_{\epsilon\to 0} \left| \chi_{B\left(0,\frac{R}{\epsilon}\right)}
    (z)g^{ij}_{x_0}(\epsilon z)|g_{x_0}(\epsilon z)|^\frac{1}{2}-
    \delta^{ij} \right| = 0
  $$
  holds true uniformly with respect to $x_0\in M$ and $z\in B(0,T)$ and 
  (\ref{limW1}) is proved.

  \noindent\textbf{(III)} \emph{There exists $t_1>0$ such that 
  $t_\epsilon\left(W_{x_0,\epsilon}\right)\geq t_1$ for any $\epsilon\in
  (0,1]$ and $x_0\in M$.}

  Let $g_{\epsilon,u}(t)=J_\epsilon(tu)$.   By Lemma~\ref{lmm-tepsilonu} 
  $(ii)$, it is enough to find $t_1>0$ such that for all $t\in[0,t_1]$ 
  $g_{\epsilon,W_{x_0,\epsilon}}'(t)>0$ for all $\epsilon\leq 1$ and for 
  all $x_0\in M$.  Then we look for a lower bound of $g_{\epsilon,W_{x_0,
  \epsilon}}'(t)$:
  $$
  \begin{aligned}
  \lefteqn{g_{\epsilon,W_{x_0,\epsilon}}'(t) = \frac{\epsilon^2t}
      {\epsilon^n} \int_M \left|\nabla W_{x_0,\epsilon}\right|^2_g 
      d\mu_g - \frac{1}{\epsilon^n} \int_M f'(tW_{x_0,\epsilon})
      W_{x_0,\epsilon}\, d\mu_g}                                    \\
    &= \frac{1}{\epsilon^n} \int_{B(0,R)} \!\big[ \epsilon^2t|\nabla 
      U_\epsilon(z)|^2_{g_{x_0}(z)} - f'(tU_\epsilon(z)-t\widetilde 
      U_\frac{R}{\epsilon}) (U_\epsilon(z)-\widetilde U_\frac{R}
      {\epsilon}) \big] |g_{x_0}(z)|^\frac{1}{2} dz                 \\
    &= \int_{B\left(0,\frac{R}{\epsilon}\right)} \big[ t|\nabla 
      U(z)|^2_{g_{x_0}(\epsilon z)} - f'(tU(z)-t\widetilde U_\frac
      {R}{\epsilon}) (U(z)-\widetilde U_\frac{R}{\epsilon}) \big] 
      |g_{x_0}(\epsilon z)|^\frac{1}{2} dz\, .
  \end{aligned}
  $$
  Using Remark~\ref{rmk-metric}, the fact that $\epsilon\leq 1$ and 
  the properties of $f$ (f1) and (f2), we obtain the following inequality:
  $$
  \begin{aligned}
  g_{\epsilon,W_{x_0,\epsilon}}'(t) > 
    &\frac{h^\frac{n}{2}t}{H} \int_{B(0,R)} |\nabla U(z)|^2 dz - 
      c_1H^\frac{n}{2} \int_{G_{t,\epsilon}} t^{p-1} \left| U(z)-
      \widetilde U_\frac{R}{\epsilon} \right|^p dz                \\
    &- c_1H^\frac{n}{2} \int_{L_{t,\epsilon}} t^{q-1} \left| U(z)-
      \widetilde U_\frac{R}{\epsilon} \right|^q dz\, ,
  \end{aligned}
  $$
  where $G_{t,\epsilon}=\left\{z\in B\left(0,\frac{R}{\epsilon}
  \right)\;|\; t\left(U(z)-\widetilde U_\frac{R}{\epsilon}\right)
  \geq 1\right\}$ and $L_{t,\epsilon}=\Big\{z\in B\left(0,\frac{R}
  {\epsilon}\right)\;|$ $\left.t\left(U(z)-\widetilde 
  U_\frac{R}{\epsilon}\right)\leq 1\right\}$.  If $t\leq 1$, the 
  following inclusions hold:
  $$
  \begin{aligned}
  G_{t,\epsilon}
    &\subset\left\{z\in B\left(0,\frac{R}{\epsilon}\right)
      \;|\; U(z)-\widetilde U_\frac{R}{\epsilon}\geq 1\right\} \\
    &\subset\left\{z\in B\left(0,\frac{R}{\epsilon}\right)
      \;|\; U(z)\geq 1\right\} \subset \left\{z\in\RR^n\;|\;
      U(z)\geq 1\right\}=G\, .
  \end{aligned}
  $$
  By these inclusions and the fact that $\left|U(z)-\widetilde U_\frac{R}
  {\epsilon}\right|\leq|U(z)|$,
  $$
  \int_{G_{t,\epsilon}} t^{p-1} \left| U(z)-\widetilde U_\frac{R}
    {\epsilon} \right|^p dz \leq \int_{G} t^{p-1} |U(z)|^p dz\, .
  $$
  Let $L=\left\{z\in\RR^n\;|\; U(z)\leq 1\right\}$. We have
  $$
  \begin{aligned}
  \int_{L_{t,\epsilon}} \hspace{-0.4cm}t^{q-1} \left| U(z)-
      \widetilde U_\frac{R}{\epsilon} \right|^q \!dz 
    &= \int_{L\cap B\left(0,\frac{R}{\epsilon}\right)} 
      \hspace{-0.7cm}t^{q-1} \left| U(z)-\widetilde U_\frac{R}
      {\epsilon} \right|^q\! dz + \int_{L_{t,\epsilon}\setminus 
      L} \hspace{-0.7cm} t^{q-1} \left| U(z)-\widetilde 
      U_\frac{R}{\epsilon} \right|^q\! dz                       \\
    &\leq \int_L t^{q-1} |U(z)|^q dz + \int_{L_{t,\epsilon}
        \setminus L} \hspace{-0.7cm} t^{p-1} \left| U(z)-
        \widetilde U_\frac{R}{\epsilon} \right|^p dz            \\
    &\leq \int_L t^{q-1} |U(z)|^q dz + \int_G  t^{p-1} |U(z)|^p 
      dz\, .
  \end{aligned}
  $$
  We conclude that
  $$
  g_{\epsilon,W_{x_0,\epsilon}}'(t) > \gamma_1 t - \gamma_2 t^{p-1} - 
    \gamma_3 t^{q-1}
  $$
  with $\gamma_1,\gamma_3$ positive constants and $\gamma_2$ nonnegative 
  constant.  This proves the existence of $t_1$.

  \noindent\textbf{(IV)} \emph{There exists $t_2>0$ such that $t_\epsilon
  \left(W_{x_0,\epsilon}\right)\leq t_2$ for any $\epsilon\in(0,1]$ and 
  $x_0\in M$.}

  If $u$ is a function in the Nehari manifold $\mathcal{N}_\epsilon$, 
  we have that $J_\epsilon(u)=\widetilde F_{\epsilon,M}(u)$, as defined 
  in (\ref{FtildeepsilonM}).   Then by property $(f1)$ $J_\epsilon(u)$ 
  is positive.    By Lemma~\ref{lmm-tepsilonu} $(ii)$, it is enough to 
  find $t_2>0$ such that for all $t\geq t_2$ $J_\epsilon(tW_{x_0,
  \epsilon})<0$ for all $\epsilon\leq 1$ and for all $x_0\in M$.  
  Then we look for an upper bound of $J_\epsilon(tW_{x_0,\epsilon})$:
  $$
  \begin{aligned}
  \lefteqn{J_\epsilon(tW_{x_0,\epsilon}) = \frac{\epsilon^2t^2}
      {2\epsilon^n} \int_M \left|\nabla W_{x_0,\epsilon}\right|^2_g 
      d\mu_g - \frac{1}{\epsilon^n} \int_M f(tW_{x_0,\epsilon})\, 
      d\mu_g}                                                      \\
    &= \frac{1}{\epsilon^n} \int_{B(0,R)} \left[ \frac{\epsilon^2
      t^2}{2}|\nabla U_\epsilon(z)|^2_{g_{x_0}(z)} - f\left(t
      U_\epsilon(z)-t\widetilde U_\frac{R}{\epsilon}\right) \right] 
      |g_{x_0}(z)|^\frac{1}{2} dz                                  \\
    &= \int_{B\left(0,\frac{R}{\epsilon}\right)} \left[ \frac{t^2}
      {2} |\nabla U(z)|^2_{g_{x_0}(\epsilon z)} - f\left(tU(z)-
      t\widetilde U_\frac{R}{\epsilon}\right) \right] |g_{x_0}
      (\epsilon z)|^\frac{1}{2} dz                                 \\
    &\leq \frac{H^\frac{n}{2}t^2}{2h}\| U\|^2_{\mathcal{D}^{1,2}
      (\RR^n)} - c_0h^\frac{n}{2} \!\!\int_{G_{t,\epsilon}} 
      \hspace{-0.4cm}t^p \left| U(z)-\widetilde U_\frac{R}
      {\epsilon} \right|^p \!\!\! dz - c_0h^\frac{n}{2} \!\!
      \int_{L_{t,\epsilon}} \hspace{-0.4cm}t^q \left| U(z)-
      \widetilde U_\frac{R}{\epsilon} \right|^q \!\!\! dz .
  \end{aligned}
  $$
  If we consider $t\geq 1$ and $\widetilde U_R=U(z)$ with $z\in\RR^n$ 
  such that $|z|=R$, there holds
  $$
  \begin{aligned}
  \int_{G_{t,\epsilon}} &t^p \left| U(z)-\widetilde 
      U_\frac{R}{\epsilon} \right|^p dz + \int_{L_{t,\epsilon}} t^q 
      \left| U(z)-\widetilde U_\frac{R}{\epsilon} \right|^q dz     \\
    \geq &\; t^p \left[ \int_{G_{1,\epsilon}} \left| U(z)-
      \widetilde U_\frac{R}{\epsilon} \right|^p dz + \int_{G_{t,
      \epsilon}\setminus G_{1,\epsilon}} \left| U(z)-\widetilde 
      U_\frac{R}{\epsilon} \right|^p dz \right.                    \\
    &+ \left. \int_{L_{1,\epsilon}} \left| U(z)-\widetilde 
      U_\frac{R}{\epsilon} \right|^q dz - \int_{L_{1,\epsilon}
      \setminus L_{t,\epsilon}} \left| U(z)-\widetilde U_\frac{R}
      {\epsilon} \right|^q dz \right]                              \\
    \geq &\; t^p \left[ \int_{G_{1,\epsilon}} \left| U(z)-
      \widetilde U_\frac{R}{\epsilon} \right|^p dz + \int_{L_{1,
      \epsilon}} \left| U(z)-\widetilde U_\frac{R}{\epsilon} 
      \right|^q dz \right]                                         \\
  \geq &\; t^p \left[ \int_{G_{1,\epsilon}\cap B(0,R)} \left| U(z)
      -\widetilde U_\frac{R}{\epsilon} \right|^p dz + \int_{L_{1,
      \epsilon}\cap B(0,R)} \left| U(z)-\widetilde U_\frac{R}
      {\epsilon} \right|^q dz \right]                              \\
    \geq &\; t^p \left[ \int_{G_{1,\epsilon}\cap B(0,R)} \left| 
      U(z)-\widetilde U_R \right|^p dz + \int_{L_{1,\epsilon}\cap 
      B(0,R)} \left| U(z)-\widetilde U_R \right|^q dz \right]      \\
     = &\; t^p \left[ \int_{G_{1,1}} \left| U(z)-\widetilde U_R 
       \right|^p dz + \int_{G_{1,\epsilon}\cap B(0,R)\setminus 
       G_{1,1}} \left| U(z)-\widetilde U_R \right|^p dz \right.     \\
     &+ \left.\int_{L_{1,1}} \left| U(z)-\widetilde U_R \right|^q dz 
       -\int_{L_{1,1}\setminus L_{1,\epsilon}} \left| 
       U(z)-\widetilde U_R \right|^q dz \right]                     \\
     \geq &\; t^p \left[ \int_{G_{1,1}} \left| U(z)-\widetilde U_R 
       \right|^p dz + \int_{L_{1,1}} \left| U(z)-\widetilde U_R 
       \right|^q dz \right].
  \end{aligned}
  $$
  So $J_\epsilon(tW_{x_0,\epsilon})\leq\gamma_4 t^2-\gamma_5 t^p$ with 
  $\gamma_4$, $\gamma_5$ positive constants and for $t$ big enough it 
  is negative.

  \noindent\textbf{(V)} \emph{The parameter $t_\epsilon\left(W_{x_0,
  \epsilon}\right)$ tends to $1$ for $\epsilon$ tending to zero uniformly 
  with respect to $x_0\in M$.}

  By the previous steps $t_\epsilon\left(W_{x_0,\epsilon}\right)\in[t_1,
  t_2]$ for any $\epsilon\in(0,1]$ and $x_0\in M$.   Let us write 
  $t_{x_0,\epsilon}=t_\epsilon\left(W_{x_0,\epsilon}\right)$.   Then 
  there exists a sequence $\epsilon_k\to 0$ for $k\to\infty$ such 
  that $t_{x_0,\epsilon_k}$ converges to $t_{x_0}^*$.  Then by step (II) 
  we have $\lim_{k\to\infty} \frac{\epsilon_k^2}{\epsilon_k^n} \int_M 
  \left| t_{x_0,\epsilon_k}\nabla W_{x_0,\epsilon_k}(x) \right|^2_g
  d\mu_g = \| t_{x_0}^*U\|^2_{\mathcal {D}^{1,2}(\RR^n)}$.   By 
  definition we have
  $$
  \begin{aligned}
  \lefteqn{\frac{1}{\epsilon_k^n} \int_M f' \left( t_{x_0,\epsilon_k}
      W_{x_0,\epsilon_k} \right) t_{x_0,\epsilon_k}W_{x_0,\epsilon_k}\, 
      d\mu_g}                                                           \\
    &= \frac{1}{\epsilon_k^n} \int_{B(0,R)} f' \big( t_{x_0,\epsilon_k}
      \big( U_{\epsilon_k}(z)-\widetilde U_\frac{R}{\epsilon_k} \big)
      \big) t_{x_0,\epsilon_k} \big( U_{\epsilon_k}(z)-\widetilde 
      U_\frac{R}{\epsilon_k} \big) |g_{x_0}(z)|^\frac{1}{2} dz          \\
    &= \int_{B(0,\frac{R}{\epsilon_k})} f' \big( t_{x_0,\epsilon_k} 
      \big( U(z)-\widetilde U_\frac{R}{\epsilon_k} \big)\big) t_{x_0,
      \epsilon_k} \big( U(z)-\widetilde U_\frac{R}{\epsilon_k} \big) 
      |g_{x_0}(\epsilon_k z)|^\frac{1}{2} dz                            \\
    &= \int_{\RR^n} \chi_{B(0,\frac{R}{\epsilon_k})}(z)f' \big( t_{x_0,
      \epsilon_k} \big( U(z)-\widetilde U_\frac{R}{\epsilon_k} \big)
      \big) t_{x_0,\epsilon_k} \big( U(z)-\widetilde U_\frac{R}
      {\epsilon_k} \big) |g_{x_0}(\epsilon_k z)|^\frac{1}{2} dz\, .
  \end{aligned}
  $$
  The integrand point-wise tends to $f'(t_{x_0}^*U(z))t_{x_0}^*U(z)$ 
  for $k$ tending to infinity and is bounded from above by a function in 
  $L^1(\RR^n)$ as follows:
  $$
  \begin{aligned}
  \lefteqn{\chi_{B(0,\frac{R}{\epsilon_k})}(z)f' \big( t_{x_0,\epsilon_k} 
      \big( U(z)-\widetilde U_\frac{R}{\epsilon_k} \big)\big) t_{x_0,
      \epsilon_k} \big( U(z)-\widetilde U_\frac{R}{\epsilon_k} \big) 
      |g_{x_0}(\epsilon_k z)|^\frac{1}{2}}                               \\
    &\leq H^\frac{n}{2}\chi_{B(0,\frac{R}{\epsilon_k})}(z)f' \big( t_2 
      \big( U(z)-\widetilde U_\frac{R}{\epsilon_k} \big)\big) t_2 \big( 
      U(z)-\widetilde U_\frac{R}{\epsilon_k} \big)                       \\
    &\leq \left\{ \begin{array}{ll}
      c_1H^\frac{n}{2}t_2^p \big( U(z)-\widetilde U_\frac{R}{\epsilon_k} 
        \big)^p &
      \mbox{if } t_2 \big( U(z)-\widetilde U_\frac{R}{\epsilon_k} \big) 
        \geq 1\mbox{ and } |z|\leq\frac{R}{\epsilon_k} \\
      c_1H^\frac{n}{2}t_2^q \big( U(z)-\widetilde U_\frac{R}{\epsilon_k} 
        \big)^q &
      \mbox{if } t_2 \big( U(z)-\widetilde U_\frac{R}{\epsilon_k} \big) 
        \leq 1\mbox{ and } |z|\leq\frac{R}{\epsilon_k} \\
      0 & \mbox{otherwise}
      \end{array} \right.                                                \\
    &\leq \left\{ \begin{array}{ll}
      \! c_1H^\frac{n}{2}t_2^p \left( U(z) \right)^p &
      \!\!\mbox{if } t_2 \big( U(z)-\widetilde U_\frac{R}{\epsilon_k} 
        \big) \geq 1,\, U(z)\geq 1\mbox{ and } |z|\leq\frac{R}
        {\epsilon_k} \\
      \! c_1H^\frac{n}{2}t_2^q \big( U(z)-\widetilde U_\frac{R}
        {\epsilon_k} \big)^q &
      \!\!\mbox{if } t_2 \big( U(z)-\widetilde U_\frac{R}{\epsilon_k} 
        \big) \geq 1,\, U(z)<1\mbox{ and } |z|\leq\frac{R}
        {\epsilon_k} \\
      \! c_1H^\frac{n}{2}t_2^p \big( U(z)-\widetilde U_\frac{R}
        {\epsilon_k} \big)^p &
      \!\!\mbox{if } t_2 \big( U(z)-\widetilde U_\frac{R}{\epsilon_k} 
        \big) <1,\, U(z)\geq 1\mbox{ and } |z|\leq\frac{R}
        {\epsilon_k} \\
      \! c_1H^\frac{n}{2}t_2^q \left( U(z) \right)^q &
      \!\!\mbox{if } t_2 \big( U(z)-\widetilde U_\frac{R}{\epsilon_k} 
        \big) <1,\, U(z)<1\mbox{ and } |z|\leq\frac{R}{\epsilon_k} \\
      \! 0 & \!\!\mbox{otherwise}
      \end{array} \right.                                                \\
    &\leq \left\{ \begin{array}{ll}
      c_1H^\frac{n}{2}t_2^p \left( U(z) \right)^p &
      \mbox{if } U(z)\geq 1 \\
      c_1H^\frac{n}{2}t_2^q \left( U(z) \right)^q &
      \mbox{if } U(z)<1
      \end{array} \right.                                                \\
    &\leq \frac{c_1H^\frac{n}{2}t_2^q}{c_0} f(U(z))\, .
  \end{aligned}
  $$
  Then by Lebesgue theorem $\lim_{k\to\infty}\frac{1}{\epsilon_k^n} 
  \int_M f' \left( t_{x_0,\epsilon_k}W_{x_0,\epsilon_k} \right) t_{x_0,
  \epsilon_k}W_{x_0,\epsilon_k}\, d\mu_g=\int_{\RR^n}f'(t_{x_0}^*U(z))
  t_{x_0}^*U(z)\, dz$.   By the fact that $U\in\mathcal{N}$ and 
  $\| t_{x_0}^*U\|^2_{\mathcal {D}^{1,2}(\RR^n)}=\int_{\RR^n}f'(t_{x_0}^*
  U(z))t_{x_0}^*U(z)\, dz$, we conclude that $t_{x_0}^*=1$.

  To prove that the convergence is uniform with respect to $x_0\in M$, 
  we show that $\lim_{\epsilon\to 0}\sup_{x\in M}|t_{x,\epsilon}-1|=0$.  
  For any $\epsilon$ there exists $x(\epsilon)\in M$ such that 
  $\sup_{x\in M}|t_{x,\epsilon}-1|=|t_{x(\epsilon),\epsilon}-1|$.  By 
  compactness there exists a sequence $\epsilon_k\to 0$ for $k\to\infty$ 
  such that $x(\epsilon_k)$ tends to $x_*\in M$.  Let us fix $\eta>0$.  
  There exists $k_0$ such that for all $k\geq k_0$ $|t_{x_*,\epsilon_k}
  -1|<\frac{\eta}{3}$.   Possibly increasing $k_0$ we also have that 
  for all $k\geq k_0$ and $h>k$ $|t_{x(\epsilon_k),\epsilon_k}-
  t_{x(\epsilon_h),\epsilon_k}|<\frac{\eta}{3}$.   Finally there exists 
  $h_0$ such that for all $h\geq h_0$ $|t_{x(\epsilon_h),\epsilon_k}-
  t_{x_*,\epsilon_k}|<\frac{\eta}{3}$.  Summing the three terms one has 
  that $|t_{x(\epsilon_k),\epsilon_k}-1|<\eta$ for all $k\geq k_0$.

  \noindent\textbf{(VI)} \emph{The limit of $\frac{1}{\epsilon^n}\int_M 
  f(t_{x_0,\epsilon}W_{x_0,\epsilon})\, d\mu_g$ is $\int_{\RR^n} f(U)\, 
  dz$}.

  Changing variables and using mean value theorem, we have
  $$
  \begin{aligned}
  \lefteqn{\frac{1}{\epsilon^n}\int_M f(t_{x_0,\epsilon}W_{x_0,
      \epsilon})\, d\mu_g = \int_{B(0,\frac{R}{\epsilon})} \left[ 
      f\big(U(z)-\widetilde U_\frac{R}{\epsilon}\big) \right.}    \\
    &+ \left.(t_{x_0,\epsilon}-1)f'\big(\Theta_{x_0,\epsilon}(z)
      \big(U(z)-\widetilde U_\frac{R}{\epsilon}\big)\big)\big(U(z)
      -\widetilde U_\frac{R}{\epsilon}\big) \right] |g_{x_0}
      (\epsilon z)|^\frac{1}{2} dz\, ,
  \end{aligned}
  $$
  where $\Theta_{x_0,\epsilon}(z)=(\theta_{x_0,\epsilon}(z)t_{x_0,
  \epsilon}+1-\theta_{x_0,\epsilon}(z))$ with a suitable $0<\theta_{x_0,
  \epsilon}(z)<1$.   We want to prove that
  \begin{equation}
  \begin{aligned}
  \int_{B(0,\frac{R}{\epsilon})} &f\big(U(z)-\widetilde U_\frac{R}
    {\epsilon}\big)|g_{x_0}(\epsilon z)|^\frac{1}{2} dz 
    \xrightarrow{\epsilon\to 0} \int_{\RR^n} f(U)\, dz              \\
  \int_{B(0,\frac{R}{\epsilon})} &\hspace{-0.6cm}(t_{x_0,\epsilon}
    -1)f'\big(\Theta_{x_0,\epsilon}(z)\big(U(z)-\widetilde U_
    \frac{R}{\epsilon}\big)\big)\big(U(z)-\widetilde U_\frac{R}
    {\epsilon}\big) |g_{x_0}(\epsilon z)|^\frac{1}{2} dz 
    \xrightarrow{\epsilon\to 0} 0
  \end{aligned}
  \label{limf}
  \end{equation}
  uniformly with respect to $x_0\in M$.
  
  It is easy to see that
  $$
  \int_{B(0,\frac{R}{\epsilon})} f\big(U(z)-\widetilde U_\frac{R}
    {\epsilon}\big) \left| |g_{x_0}(\epsilon z)|^\frac{1}{2}-1 \right| dz
    \xrightarrow{\epsilon\to 0} 0
  $$
  uniformly with respect to $x_0\in M$.   The function $\chi_{B\left(0,
  \frac{R}{\epsilon}\right)}(z)f\big(U(z)-\widetilde U_\frac{R}{\epsilon}
  \big)$ tends pointwise to $f(U(z))$ for any $z\in\RR^n$.  Moreover
  $$
  \begin{aligned}
  \chi_{B\left(0,\frac{R}{\epsilon}\right)}(z)f\big(U(z)
      -\widetilde U_\frac{R}{\epsilon}\big)
    &\leq \left\{ \begin{array}{ll} 
      \frac{c_1}{\mu} \big(U(z)-\widetilde U_\frac{R}{\epsilon}\big)^p &
      \mbox{if } U(z)-\widetilde U_\frac{R}{\epsilon} \geq 1,\;
        |z|\leq\frac{R}{\epsilon} \\
      \frac{c_1}{\mu} \big(U(z)-\widetilde U_\frac{R}{\epsilon}\big)^q &
      \mbox{if } U(z)-\widetilde U_\frac{R}{\epsilon} \leq 1,\;
        |z|\leq\frac{R}{\epsilon} \\
      0 & \mbox{otherwise}
      \end{array} \right.                                               \\
    &\leq \left\{ \begin{array}{ll} 
      \frac{c_1}{\mu} \big(U(z)-\widetilde U_\frac{R}{\epsilon}\big)^p &
      \mbox{if } U(z)\geq 1,\; |z|\leq\frac{R}{\epsilon}\\
      \frac{c_1}{\mu} \big(U(z)-\widetilde U_\frac{R}{\epsilon}\big)^q &
      \mbox{if } U(z)<1,\; |z|\leq\frac{R}{\epsilon}\\
      0 & \mbox{otherwise}
      \end{array} \right.                                                \\
    &\leq \left\{ \begin{array}{ll} 
      \frac{c_1}{\mu} (U(z))^p & \mbox{if } U(z)\geq 1 \\
      \frac{c_1}{\mu} (U(z))^q & \mbox{if } U(z)\leq 1 
      \end{array} \right.                                                \\
    &\leq \frac{c_1}{c_0\mu} f(U(z))
  \end{aligned}
  $$
  and by Lebesgue theorem we obtain the first limit in (\ref{limf}).   
  The function of $t$ $f'(tu)u$ is increasing in $t$, since its 
  derivative is $f''(tu)u^2>0$.  Then we have
  $$
  \begin{aligned}
  \int_{B(0,\frac{R}{\epsilon})}
    &f'\big(\Theta_{x_0,\epsilon}(z)\big(U(z)-\widetilde U_
      \frac{R}{\epsilon}\big)\big)\big(U(z)-\widetilde U_
      \frac{R}{\epsilon}\big) |g_{x_0}(\epsilon z)|^\frac{1}{2} dz \\
    &< H^\frac{n}{2}\int_{B(0,\frac{R}{\epsilon})} f'\big((t_2+1)
      \big(U(z)-\widetilde U_\frac{R}{\epsilon}\big)\big)\big(U(z)
      -\widetilde U_\frac{R}{\epsilon}\big) dz\, .
  \end{aligned}
  $$
  By the usual standard inequalities, the previous integral is 
  bounded from above by $\frac{c_1H^\frac{n}{2}}{c_0(t_2+1)}
  \int_{\RR^n}f((t_2+1)U(z))dz$ and the second limit in (\ref{limf}) 
  is proved, because of (V).

  \noindent\textbf{(VII)} \emph{Conclusion.}

  By (II), (V) and (VI) we obtain that $J_\epsilon(\phi_\epsilon(x_0))$ 
  tends to $J(U)=m(J)$ for $\epsilon$ tending to zero uniformly with 
  respect to $x_0$. This completes the proof.
\end{proof}

\begin{remark}
  By the previous proposition, in particular we know that, given 
  $\delta>0$, for any positive $\epsilon$ sufficiently small 
  $\Sigma_{\epsilon,\delta}$ is not empty.
\end{remark}

\section{The function $\beta$}

Given a function $u\in L^p(M)$, $u\not\equiv 0$, it is possible to 
define its centre of mass $\beta(u)\in\RR^N$ by
\begin{equation}
\beta(u) = \frac{\int_Mx\Phi(u)\, d\mu_g}{\int_M\Phi(u)\, d\mu_g}\, ,
\label{beta}
\end{equation}
where
\begin{equation}
\Phi(u) = \frac{1}{2}f'(u)u-f(u)\, .
\label{Phi}
\end{equation}
By the properties of $f$, $\Phi(s)>0$ for all $s\neq 0$.   To prove 
that $\beta:\Sigma_{\epsilon,\delta}\to M_{r(M)}$ (see 
Section~\ref{sec:ideas} and Definition~\ref{def-r(M)}), we use the 
fact that the functions in $\Sigma_{\epsilon,\delta}$ concentrate 
for $\epsilon$ and $\delta$ tending to zero.   

First of all we find a positive inferior bound for the functional 
$J_\epsilon$ on the Nehari manifold.  Let us denote
\begin{equation}
m_\epsilon = \inf_{u\in\mathcal{N}_\epsilon} J_\epsilon(u)\, .
\label{mepsilon}
\end{equation}

It is easy to see that
$$
\inf_{u\in\mathcal{N}_\epsilon}\| u\|_{H^1_2(M)}>0
$$
(the proof is analogous to Lemma 3.2 of \cite{BM}) and, since the 
manifold $M$ is compact, that the infimum $m_\epsilon$ is achieved.

\begin{lemma}
\label{lmm-mepsilon>alpha}
  There exist positive constants $\alpha$ and $\epsilon_0$ such that 
  for any $0<\epsilon<\epsilon_0$ the inequality $m_\epsilon\geq\alpha$ 
  holds.
\end{lemma}

To prove this lemma we need the following technical lemma (for the 
proof see the Appendix).

\begin{lemma}
\label{lmm-estensione}
  For any $r\in(0,r(M))$, there exist constants $k_1,\,k_2,\,k_3,\,k_4
  >0$ such that for any $u\in H^1_2(M)$ there exists $v\in
  \mathcal{D}^{1,2}(M_r)$ such that $v|_M\equiv u$ and
  \begin{eqnarray}
  \|v\|^2_{\mathcal{D}^{1,2}(M_r)} 
    & \leq & k_1 \int_M |\nabla u|^2_g d\mu_g\, , \label{k1} \\
  \int_{M_r}f(v(z))\, dz 
    & \geq & k_2 \int_M f(u(x))\, d\mu_g\, , \label{k2}      \\
  \int_{M_r}f(v(z))\, dz 
    & \leq & k_3 \int_M f(u(x))\, d\mu_g\, , \label{k3}      \\
  \|v\|^2_{L^2(M_r)}
    & \geq & k_4 \|u\|^2_{L^2(M)}\, . \label{k4}
  \end{eqnarray}
\end{lemma}

\begin{proof}[Proof of Lemma~\ref{lmm-mepsilon>alpha}]
  By definition $m_\epsilon$ is the infimum of $J_\epsilon(u)$ on the 
  Nehari manifold $\mathcal{N}_\epsilon$.  If $u\in\mathcal{N}_\epsilon$ 
  we have
  $$
  J_\epsilon(u) \geq \left(\frac{1}{2}-\frac{1}{\mu} \right) 
    \frac{\epsilon^2}{\epsilon^n} \int_M |\nabla u|^2_g d\mu_g\, .
  $$
  Rescaling $u$, it is easy to see that $m_\epsilon$ is greater than or 
  equal to the infimum of the functional $\left(\frac{1}{2}-\frac{1}{\mu} 
  \right)\frac{\epsilon^2}{\epsilon^n}t_\epsilon^2 \int_M |\nabla w|^2_g 
  d\mu_g$ on the set of the functions $w\in H^1_2(M)$ such that 
  $\frac{1}{\epsilon^n} \int_M f(w)\, d\mu_g=1$ and where $t_\epsilon=
  t_\epsilon(w)$ is as in (ii), Lemma~\ref{lmm-tepsilonu}.  First of all, 
  we check that there exists a constant $\tilde\alpha>0$ and for such 
  functions $w$ it holds
  $$
  \frac{\epsilon^2}{\epsilon^n} \int_M |\nabla w|^2_g d\mu_g \geq 
    \tilde\alpha\, .
  $$
  By Lemma \ref{lmm-estensione}, for any function $w$ there exists a 
  function $v\in\mathcal{D}^{1,2}(M_r)$ such that (\ref{k1}) and 
  (\ref{k2}) hold.   We consider $\tilde v\in\mathcal{D}^{1,2}(\RR^N)$, 
  defined as $\tilde v(y)=v(y)$ for all $y\in M_r$ and $\tilde v(y)=0$ 
  for all $y\in\RR^N\setminus M_r$.  We can now consider the following 
  rescalement $V(y)=\tilde v(\epsilon^\sigma y)$ with $\sigma=
  \frac{2n-(n-2)p}{2N-(N-2)p}$.  In case the denominator is equal to $0$, 
  we can choose a bigger $N$.  We have:
  $$
  \| V\|^2_{\mathcal{D}^{1,2}(\RR^N)} = 
    \frac{\epsilon^{2\sigma}}{\epsilon^{N\sigma}} \| v\|^2
    _{\mathcal{D}^{1,2}(M_r)} \ \mbox{and}\ \int_{\RR^N} f(V(y))\, dy 
    = \frac{1}{\epsilon^{N\sigma}} \int_{M_r} f(v(y))\, dy\, .
  $$
  By these equalities, (\ref{k1}) and (\ref{k2}), we have
  \begin{equation}
  \begin{aligned}
  \frac{\epsilon^2}{\epsilon^n} \int_M |\nabla w|^2_g d\mu_g 
    &= \frac{\frac{\epsilon^2}{\epsilon^n} \int_M |\nabla w|^2_g 
      d\mu_g}{\left(\frac{1}{\epsilon^n}\int_M f(w)\, d\mu_g\right)^
      \frac{2}{p}} \geq \frac{k_2^\frac{2}{p}}{k_1} 
      \frac{\frac{\epsilon^2}{\epsilon^n} \| v\|^2_{\mathcal{D}^{1,2}
      (M_r)}}{\left(\frac{1}{\epsilon^n}\int_{M_r} f(v)\, dy\right)^
      \frac{2}{p}}                                                   \\
    &= \frac{k_2^\frac{2}{p}}{k_1} \frac{\frac{\epsilon^{(N-2)
      \sigma}}{\epsilon^{n-2}} \| V\|^2_{\mathcal{D}^{1,2}(\RR^N)}}
      {\left(\frac{\epsilon^{N\sigma}}{\epsilon^n}\int_{\RR^N} f(V)
      \, dy\right)^\frac{2}{p}} = \frac{k_2^\frac{2}{p}}{k_1} 
      \frac{\| V\|^2_{\mathcal{D}^{1,2}(\RR^N)}}{\left(\int_{\RR^N} 
      f(V)\, dy\right)^\frac{2}{p}}\, .
  \end{aligned}
  \label{dis1}
  \end{equation}
  We show now that for $\epsilon$ sufficiently small we have 
  $\int_{\RR^N} f(V)\, dy<1$.  In fact, by (\ref{k3}) there holds
  $$
  \int_{\RR^N} f(V)\, dy = \frac{1}{\epsilon^{N\sigma}} \int_{M_r} 
    f(v(y))\, dy \leq \frac{k_3}{\epsilon^{N\sigma}} \int_{M} f(w)
    \, d\mu_g = \frac{k_3\epsilon^n}{\epsilon^{N\sigma}}\, .
  $$
  By definition of $\sigma$ $\lim_{N\to\infty} N\sigma=\frac{2n-
  (n-2)p}{2-p}<0$ and so there exists $N$ sufficiently big such that 
  $n-N\sigma>0$.
  
  Since $\int_{\RR^N} f(tV(y))\, dy$ is an increasing function of 
  $t$ for positive $t$, there exists $t_*>1$ such that $\int_{M_r} f(t_*V(y))
  \, dy=1$.  Let $V_*(y)=t_*V(y)$ for any $y\in\RR^N$.   With the 
  usual computation we obtain
  $$
  \begin{aligned}
  \int_{\RR^N} &f(V(y))\, dy = \int_{\RR^N} f \left( 
      \frac{1}{t_*}V^*(y) \right) dy                          \\
    < &\frac{c_1}{\mu} \left( \int_{\{ y\in\RR^N\;|\;|V_*(y)|
      \geq t_*\}} \frac{1}{t_*^p}|V_*(y)|^p dy + \int_{\{ y\in
      \RR^N\;|\;|V_*(y)|\leq t_*\}} \frac{1}{t_*^q}|V_*(y)|^q 
      dy \right)                                              \\
    \leq &\frac{c_1}{\mu} \left( \int_{\{ y\in\RR^N\;|\;
      |V_*(y)|\geq 1\}} \frac{1}{t_*^p}|V_*(y)|^p dy + \int_{
      \{ y\in\RR^N\;|\;|V_*(y)|\leq 1\}} \frac{1}{t_*^q}|V_*
      (y)|^q dy \right)                                       \\
    \leq &\frac{c_1}{c_0\mu t_*^p} \int_{\RR^N} f(V_*(y))\, dy 
      = \frac{c_1}{c_0\mu t_*^p}\, .
  \end{aligned}
  $$
  Concluding we have that the last term in (\ref{dis1}) is equal to
  $$
  \frac{k_2^\frac{2}{p}}{k_1} \frac{\frac{1}{t_*^2}\| V_*\|^2_
    {\mathcal{D}^{1,2}(\RR^N)}}{\left(\int_{\RR^N} f \left( 
    \frac{1}{t_*}V_* \right) dy\right)^\frac{2}{p}} \geq 
    \frac{k_2^\frac{2}{p}}{k_1} \left( \frac{c_0\mu}{c_1} 
    \right)^\frac{2}{p} \| V_*\|^2_{\mathcal{D}^{1,2}(\RR^N)}
  $$
  which is bounded from below because (see \cite{BL})
  $$
  \inf_{V\in\mathcal{D}^{1,2}(\RR^N)\atop\int_{\RR^N} f(V)\, dy=1} 
    \| V\|^2_{\mathcal{D}^{1,2}(\RR^N)}=\hat\alpha>0\, .
  $$

  We still have to show that $t_\epsilon$ is bounded from below by 
  a positive constant.  By the properties (f1) and (f2) we have
  $$
  \begin{aligned}
  \frac{1}{\epsilon^n} \int_M f'(t_\epsilon w)t_\epsilon w\, d\mu_g 
    &< \frac{c_1}{\epsilon^n} \left[ \int_{\{x\in M\;|\;|t_\epsilon 
      w(x)|\geq 1\}} \hspace{-1.5cm}|t_\epsilon w(x)|^p d\mu_g + 
      \int_{\{x\in M\;|\;|t_\epsilon w(x)|\leq 1\}} \hspace{-1.5cm}
      |t_\epsilon w(x)|^q d\mu_g \right]                            \\
    &\leq \frac{c_1}{\epsilon^n} \left[ \int_{\{x\in M\;|\;|w(x)|
      \geq 1\}} \hspace{-1.5cm}|t_\epsilon w(x)|^p d\mu_g + 
      \int_{\{x\in M\;|\;|w(x)|\leq 1\}} \hspace{-1.5cm}|t_\epsilon 
      w(x)|^q d\mu_g \right]                                        \\
    &\leq \frac{c_1t_\epsilon^p}{c_0\epsilon^n} \int_M f(w(x))\, 
      d\mu_g = \frac{c_1t_\epsilon^p}{c_0}\, ,
  \end{aligned}
  $$
  where the last equality is due the property of the functions $w$.  
  Since $t_\epsilon w\in\mathcal{N}_\epsilon$, $\frac{1}{\epsilon^n} 
  \int_M f'(t_\epsilon w)t_\epsilon w\, d\mu_g=\frac{\epsilon^2
  t_\epsilon^2}{\epsilon^n} \int_M |\nabla w|^2_g d\mu_g$ and by the 
  previous inequalities we have
  $$
  t_\epsilon^{p-2} \geq \frac{c_0}{c_1} \frac{\epsilon^2}{\epsilon^n} 
    \int_M |\nabla w|^2_g d\mu_g \geq \frac{c_0}{c_1}\widetilde\alpha
  $$
  and this completes the proof.
\end{proof}


In the following lemma for every function $u\in\mathcal{N}_\epsilon$ 
it is stated the existence of a point in the manifold where $u$ in 
some sense concentrates.

\begin{lemma}
\label{lmm-gamma}
  Let $\mathcal{C}$ be an atlas for $M$ with open cover given by 
  $B_g(x_i,R)$, $i=1,\dots ,A$, and partition of unity 
  $\{\psi_i\}_{i=1\dots A}$.  There exists a constant $\gamma>0$ 
  such that for any $0<\epsilon<\epsilon_0$, where $\epsilon_0$ 
  is defined in Lemma \ref{lmm-mepsilon>alpha}, if $u\in
  \mathcal{N}_\epsilon$ there exists $i=i(u)$ such that
  \begin{equation}
  \begin{aligned}
  \frac{1}{\epsilon^n} \int_{B_g\left(x_i,\frac{R}{2}\right)} 
    \left[ \frac{1}{2}f'(u)u-f(u) \right] d\mu_g &\geq \gamma\, , \\
  \frac{\epsilon^2}{2\epsilon^n} \int_{B_g\left(x_i,\frac{R}{2}
    \right)} |\nabla u|^2_g\, d\mu_g - \frac{1}{\epsilon^n} 
    \int_{B_g\left(x_i,\frac{R}{2}\right)} f(u)\, d\mu_g &\geq 
    \gamma\, .
  \end{aligned}
  \label{gamma}
  \end{equation}
\end{lemma}

\begin{proof}
  Let $u$ be in $\mathcal{N}_\epsilon$.    We assume that 
  $\widetilde{\mathcal{C}}=\left\{ B_g\left(x_i,\frac{R}{2}
  \right)\right\}_{i=1,\dots,A}$ is still an open cover 
  (otherwise we complete $\mathcal{C}$).   Let $\{\tilde\psi_i
  \}_{i=1\dots A}$ be a partition of unity subordinate to the 
  atlas $\widetilde{\mathcal{C}}$.   If $\widetilde F_{\epsilon,
  M}(u)$ is as in (\ref{FtildeepsilonM}), it is possible 
  to write
  $$
  \begin{aligned}
  \lefteqn{J_\epsilon(u) = \left( \widetilde F_{\epsilon,M}(u) 
      \right)^\frac{1}{2} \left( J_\epsilon(u) \right)^\frac{1}{2}} \\
    &= \left( \frac{1}{\epsilon^n} \sum_{i=1}^A \int_{B_g\left(x_i,
      \frac{R}{2}\right)} \tilde\psi_i(x) \left[ \frac{1}{2}f'
      (u(x))u(x)-f(u(x)) \right] d\mu_g \right)^\frac{1}{2} \left( 
      J_\epsilon(u) \right)^\frac{1}{2}                             \\
    &\leq \sqrt{A}\max_{1\leq i\leq A} \left( \widetilde 
      F_{\epsilon,B_g\left(x_i,\frac{R}{2}\right)}(u) \right)^
      \frac{1}{2} \left( J_\epsilon(u) \right)^\frac{1}{2}          
  \end{aligned}
  $$
  By this inequality and Lemma \ref{lmm-mepsilon>alpha} we conclude that
  $$
  \max_{1\leq i\leq A} \widetilde F_{\epsilon,B_g\left(x_i,\frac{R}{2}
    \right)}(u)\geq \frac{1}{A} J_\epsilon(u) \geq \frac{\alpha}{A}\, .
  $$
  The second equation in (\ref{gamma}) is proved analogously.
\end{proof}

In the following proposition the concentration property is better 
specified.

\begin{proposition}
\label{prp-x0}
  For any $\eta\in(0,1)$ there exists $\delta_0<m(J)$ such that, for 
  any $\delta\in(0,\delta_0)$ there exists $\epsilon_0=\epsilon_0
  (\delta)>0$ and for any $\epsilon\in(0,\epsilon_0)$ with every 
  function $u\in\Sigma_{\epsilon,\delta}$ it is associated a point 
  $x_0=x_0(u)$ in $M$ with the property
  $$
  \widetilde F_{\epsilon,B_g\left(x_0,\frac{r(M)}{2}\right)}(u) > 
    \eta\, m(J)\, .
  $$
\end{proposition}

The proof of this proposition needs the following lemmas.  The first 
lemma we need is the splitting lemma proved in \cite{BM} (Lemma 4.1):

\begin{lemma}
\label{lmm-splitting}
  Let $\{v_k\}_{k\in\NN}\subset\mathcal{N}$ be a sequence such that:
  $$
  \begin{array}{ll}
    J(v_k)\to m(J) & \mbox{as }k\to\infty\, , \\
    J'(v_k)\to 0 \mbox{ in }\mathcal{D}^{1,2}(\RR^n) & 
    \mbox{as }k\to\infty\, .
    \end{array}
  $$
  Then
  \begin{itemize}
  \item either $v_k$ converges strongly in $\mathcal{D}^{1,2}(\RR^n)$ 
    to a ground state solution of (\ref{eqRn})
  \item or there exist a sequence of points $\{y_k\}_{k\in\NN}
    \subset\RR^n$ with $|y_k|\to\infty$ as $k\to\infty$, a ground 
    state solution $U$ of (\ref{eqRn}) and a sequence of functions 
    $\{v^0_k\}_{k\in\NN}$ such that, up to a subsequence:\\
    (i) $v_k(z)=v^0_k(z)+U(z-y_k)$ for all $z\in\RR^n$;\\
    (ii) $v^0_k\to 0$ as $k\to\infty$ in $\mathcal{D}^{1,2}(\RR^n)$.
  \end{itemize}
\end{lemma}

\begin{lemma}
\label{limSigma}
  Let $\epsilon_k$ and $\delta_k$ be two positive sequences 
  tending to zero for $k$ tending to infinity.  For any $k\in\NN$ 
  let $u_k$ be a function in $\Sigma_{\epsilon_k,\delta_k}$ such that 
  for any $u\in H^1_2(M)$
  $$
  |J_{\epsilon_k}'(u_k)(u)| = o \left( \frac{\epsilon_k}{\epsilon_k^
    \frac{n}{2}} \| u\|_{H^1_2(M)} \right).
  $$
  There exist a sequence $\{x_k\}_{k\in\NN}$ of points in $M$ and a 
  sequence of functions $w_k$ on $\RR^n$, defined as
  \begin{equation}
  w_k(z) = u_k(\exp_{x_k}(\epsilon_k z)) \chi_\frac{R}{\epsilon_k}
    (|z|)\, ,
  \label{w_k}
  \end{equation}
  such that the following properties hold:
  \begin{itemize}
  \item[(i)] There exists $w\in\mathcal{D}^{1,2}(\RR^n)$ such that, up 
    to a subsequence, $w_k$ tends to $w$ weakly in $\mathcal{D}^{1,2}
    (\RR^n)$ and strongly in $L^p_{loc}(\RR^n)$.
  \item[(ii)] The function $w$ is a weak solution of $-\Delta w=f'(w)$ 
    on $\RR^n$.
  \item[(iii)] The function $w$ is a ground state solution.
  \item[(iv)]  The following equality holds
    $$
    \lim_{k\to\infty} J_{\epsilon_k}(u_k) = m(J)\, .
    $$
  \end{itemize}
\end{lemma}

\begin{proof}
  To get started we consider $x_k$ to be the points in $M$ such that 
  $u_k$ has the property (\ref{gamma}).  We will be more precise in 
  point $(iii)$.\\
  $(i)$ It is sufficient to prove that the sequence $w_k$ is bounded in 
  $\mathcal{D}^{1,2}(\RR^n)$.   We write:
  $$
  \begin{aligned}
  \|w_k\|^2_{\mathcal{D}^{1,2}(\RR^n)} = &\int_{B\left(0,
      \frac{R}{\epsilon_k}\right)} |\nabla w_k(z)|^2 dz       \\
    \leq& 2\int_{B\left(0,\frac{R}{\epsilon_k}\right)} 
      |\nabla[u_k(\exp_{x_k}(\epsilon_kz))]|^2 \left[ 
      \chi_\frac{R}{\epsilon_k}(|z|) \right]^2 dz             \\
    &+ 2\int_{B\left(0,\frac{R}{\epsilon_k}\right)} \left[
      \chi_\frac{R}{\epsilon_k}'(|z|)\right]^2 [u_k(\exp_{x_k}
      (\epsilon_kz))]^2 dz = I_1+I_2\, .
  \end{aligned}
  $$
  We consider the following inequality:
  \begin{equation}
  \begin{aligned}
  \frac{\epsilon_k^2}{\epsilon_k^n} \int_M |\nabla u_k|_g^2 d\mu_g 
    &\geq \frac{\epsilon_k^2}{\epsilon_k^n} \int_{B_g(x_k,R)} 
      |\nabla u_k|_g^2 d\mu_g                                     \\
    &= \frac{\epsilon_k^2}{\epsilon_k^n} \int_{B(0,R)} |\nabla 
      u_k(\exp_{x_k}(z))|^2_{g_{x_k}(z)} |g_{x_k}(z)|^\frac{1}{2} 
      dz                                                          \\
    &= \int_{B\left(0,\frac{R}{\epsilon_k}\right)} |\nabla 
      u_k(\exp_{x_k}(\epsilon_k z))|^2_{g_{x_k}(\epsilon_k z)} 
      |g_{x_k}(\epsilon_k z)|^\frac{1}{2} dz                      \\
    &\geq \frac{h^\frac{n}{2}}{H} \int_{B\left(0,
      \frac{R}{\epsilon_k}\right)} |\nabla u_k(\exp_{x_k}
      (\epsilon_k z))|^2 dz \geq \frac{h^\frac{n}{2}}{2H}\, I_1\, .
  \end{aligned}
  \label{dis-I1}
  \end{equation}
  Moreover the following inequality holds
  \begin{equation}
  \begin{aligned}
  I_2 &\leq \frac{2\chi_0^2\epsilon_k^2}{R^2} \int_{B\left(0,
      \frac{R}{\epsilon_k}\right)} [u_k(\exp_{x_k}(\epsilon_k
      z))]^2 dz                                               \\
    &= \frac{2\chi_0^2\epsilon_k^2}{R^2\epsilon_k^n} 
      \int_{B(0,R)} [u_k(\exp_{x_k}(z))]^2 dz                 \\
    &\leq \frac{2\chi_0^2\epsilon_k^2}{h^\frac{n}{2}R^2
      \epsilon_k^n} \int_{B_g(x_k,R)} (u_k(x))^2 d\mu_g\, .
  \end{aligned}
  \label{dis-I2}
  \end{equation}
  By (\ref{dis-I1}) and (\ref{dis-I2}), we have that the sum $I_1+I_2$ 
  is bounded by a constant times $\frac{\epsilon_k^2}{\epsilon_k^n}\| 
  u_k\|^2_{H^1_2(M)}$.   We show then that this quantity must be 
  bounded.  Since $u_k\in\Sigma_{\epsilon_k,\delta_k}$ and
  $$
  J_{\epsilon_k}(u_k) \geq \left( \frac{1}{2}-\frac{1}{\mu} \right) 
    \frac{\epsilon_k^2}{\epsilon_k^n} \int_M |\nabla u_k|_g^2 d\mu_g\, ,
  $$
  the right hand side of the preceding inequality must be bounded.  We 
  still have to check that $\frac{\epsilon_k^2}{\epsilon_k^n}\| 
  u_k\|^2_{L^2(M)}$ is bounded too.  In fact, by (\ref{k4}) in 
  Lemma~\ref{lmm-estensione} we have a sequence $v_k$ of functions 
  in $\mathcal{D}^{1,2}(M_r)$ and
  $$
  \frac{\epsilon_k^2}{\epsilon_k^n}\| u_k\|^2_{L^2(M)} \leq 
    \frac{\epsilon_k^2}{k_4\epsilon_k^n}\| v_k\|^2_{L^2(M_r)} 
    \leq \frac{C\epsilon_k^2}{k_4\epsilon_k^n}\| v_k\|^2_
    {\mathcal{D}^{1,2}(M_r)} \leq \frac{Ck_1\epsilon_k^2}
    {k_4\epsilon_k^n} \int_M |\nabla u_k|^2_g d\mu_g\, ,
  $$
  where $C$ is the constant in the Poincar\'e inequality and we have 
  used (\ref{k1}) in the last inequality.

  \noindent $(ii)$ First of all we prove that for any $\xi\in 
  C_0^\infty(\RR^n)$ $J'(w_k)(\xi)$ tends to zero for $k$ tending to 
  infinity:
  $$
  \begin{aligned}
  \lefteqn{J'(w_k)(\xi) = \int_{\RR^n} \nabla w_k(z)\cdot\nabla
      \xi(z)\, dz - \int_{\RR^n} f'(w_k(z)) \xi(z)\, dz}           \\
    &= \!\!\!\int_{\RR^n} \!\!\left[ \nabla \!\big[ u_k(\exp_{x_k}
      (\epsilon_k z)) \chi_\frac{R}{\epsilon_k}(|z|) \big] \!
      \cdot\!\nabla\xi(z) - f'\!\big( u_k(\exp_{x_k}(\epsilon_k 
      z))\chi_\frac{R}{\epsilon_k}(|z|) \big) \xi(z) \right]\!\! 
      dz                                                           \\
    &= \!\!\!\int_{\RR^n} \left[ \nabla \left[ u_k(\exp_{x_k}
      (\epsilon_k z)) \right] \cdot \nabla\xi(z) - f' \left( u_k
      (\exp_{x_k}(\epsilon_k z)) \right) \xi(z) \right] dz\, ,
  \end{aligned}
  $$
  where in the last equality we have used the fact that for $k$ 
  sufficiently large for any $z$ in the support of $\xi$ $\chi_
  \frac{R}{\epsilon_k}(|z|)=1$.  Now we define the function 
  $\tilde\xi_k$ in $H^1_2(M)$ as follows:
  $$
  \tilde\xi_k(x) = \left\{ \begin{array}{ll}
    \xi \left( \frac{\exp_{x_k}^{-1}(x)}{\epsilon_k} \right) & 
    \forall x\in B_g(x_k,R)\, ,                              \\
    0 & \mbox{otherwise.}
    \end{array} \right.
  $$
  Then we want to write
  $$
  J'(w_k)(\xi) = \frac{\epsilon_k^2}{\epsilon_k^n} \int_M g_{x_k} 
    \left( \nabla u_k(x),\nabla\tilde\xi_k(x) \right) d\mu_g - 
    \frac{1}{\epsilon_k^n} \int_M f'(u_k(x))\tilde\xi_k(x)\, d\mu_g 
    + E_k\, ,
  $$
  where $E_k$ is an error.   By hypothesis
  $$
  \begin{aligned}
  \lefteqn{\left| \int_M \left[ \frac{\epsilon_k^2}{\epsilon_k^n} 
      g_{x_k} \left( \nabla u_k(x),\nabla\tilde\xi_k(x) \right) - 
      \frac{1}{\epsilon_k^n} f'(u_k(x))\tilde\xi_k(x) \right] 
      d\mu_g \right|}                                             \\
    &= \left| J_{\epsilon_k}'(u_k)(\tilde\xi_k) \right| = o \left( 
      \frac{\epsilon_k}{\epsilon_k^\frac{n}{2}} \|\tilde\xi
      \|_{H^1_2(M)} \right) = o(\|\xi\|_{H^1_2(\RR^n)}).
  \end{aligned}
  $$
  Now we have to check the error:
  $$
  \begin{aligned}
  |E_k| = & \left| \int_{\RR^n} \left[ \nabla \left[ u_k(\exp_{x_k}
      (\epsilon_k z)) \right] \cdot \nabla\xi(z) - f' \left( u_k
      (\exp_{x_k}(\epsilon_k z)) \right) \xi(z) \right] dz \right. \\
    &- \left. \frac{\epsilon_k^2}{\epsilon_k^n} \int_M g_{x_k} 
      \left( \nabla u_k(x),\nabla\tilde\xi_k(x) \right) d\mu_g - 
      \frac{1}{\epsilon_k^n} \int_M f'(u_k(x))\tilde\xi_k(x)\, 
      d\mu_g \right|                                               \\
    \leq & \left| \int_{\RR^n} \nabla \left[ u_k(\exp_{x_k}
      (\epsilon_k z)) \right] \cdot \nabla\xi(z)\, dz - 
      \frac{\epsilon_k^2}{\epsilon_k^n} \int_M g_{x_k} \left( 
      \nabla u_k(x),\nabla\tilde\xi_k(x) \right) d\mu_g \right|    \\
    &+ \left| \int_{\RR^n} f' \left( u_k(\exp_{x_k}(\epsilon_k z)) 
      \right) \xi(z)\, dz - \frac{1}{\epsilon_k^n} \int_M f'(u_k
      (x))\tilde\xi_k(x)\, d\mu_g \right|                          \\
    = & |E_{1,k}|+|E_{2,k}|\, .
  \end{aligned}
  $$
  For the first term we have
  $$
  |E_{1,k}| \leq \int_\Xi \left| (\delta^{ij}-g^{ij}_{x_k}
    (\epsilon_k z)|g_{x_k}(\epsilon_k z)|^\frac{1}{2}) 
    \frac{\partial [u_k(\exp_{x_k}(\epsilon_k z))]}{\partial z_i}
    \frac{\partial\xi(z)}{\partial z_j} \right| dz\, ,
  $$
  where $\Xi$ denotes the compact support of $\xi$.   The limit 
  $$
  \lim_{k\to\infty}|\delta^{ij}-g^{ij}_{x_k}(\epsilon_k z)
    |g_{x_k}(\epsilon_k z)|^\frac{1}{2}|=0
  $$
  is uniform with respect to $z\in\Xi$.  Since
  $$
  \int_\Xi \left| \frac{\partial [u_k(\exp_{x_k}(\epsilon_k z))]}
    {\partial z_i}\frac{\partial\xi(z)}{\partial z_j} \right| dz 
    \leq \| u_k(\exp_{x_k}(\epsilon_k z))\|_{\mathcal{D}^{1,2}(\Xi)}
    \|\xi\|_{\mathcal{D}^{1,2}(\RR^n)}
  $$
  and for $k$ sufficiently large
  $$
  \begin{aligned}
  \int_\Xi \left| \nabla [u_k(\exp_{x_k}(\epsilon_k z))] \right|^2 dz 
    &\leq \frac{H}{h^\frac{n}{2}} \frac{\epsilon_k^2}{\epsilon_k^n} 
      \int_M |\nabla u_k|_g^2 d\mu_g                                 \\
    &\leq \frac{2\mu H}{(\mu-2)h^\frac{n}{2}} J_{\epsilon_k}(u_k) 
      \leq \frac{4\mu Hm(J)}{(\mu-2)h^\frac{n}{2}}\, ,
  \end{aligned}
  $$
  we conclude that $|E_{1,k}|$ tends to zero.   For the second term 
  we have
  $$
  |E_{2,k}| = \left| \int_\Xi \left( 1-|g_{x_k}(\epsilon_k z)|^
    \frac{1}{2} \right) f' \left( u_k(\exp_{x_k}(\epsilon_k z)) 
    \right) \xi(z)\, dz \right|.
  $$
  As before, $\lim_{k\to\infty}|g_{x_k}(\epsilon_k z)|^\frac{1}{2}$ is 
  $1$ uniformly with respect to $z\in\Xi$ and
  $$
  \begin{aligned}
  \int_\Xi &\left| f' \left( u_k(\exp_{x_k}(\epsilon_k z)) \right) 
      \xi(z)\right| dz                                            \\
    \leq &\left( \int_{\{z\in\Xi\;|\;|u_k(\exp_{x_k}(\epsilon_k 
      z))|\geq 1\}} \left| f' \left( u_k(\exp_{x_k}(\epsilon_k z)) 
      \right)\right|^\frac{p}{p-1} dz \right)^\frac{p-1}{p} 
      \|\xi\|_{L^p(\RR^n)}                                        \\
    &+ \left( \int_{\{z\in\Xi\;|\;|u_k(\exp_{x_k}(\epsilon_k 
      z))|\leq 1\}} \left| f' \left( u_k(\exp_{x_k}(\epsilon_k z)) 
      \right)\right|^\frac{q}{q-1} dz \right)^\frac{q-1}{q} 
      \|\xi\|_{L^q(\RR^n)}\, .
  \end{aligned}
  $$
  It is easy to see that there exists a positive constant $C$ such that 
  the right side is bounded from above by
  $$
  \begin{aligned}
  C &\left[ \left( \frac{1}{\epsilon_k^n} \int_M f'(u_k)u_k\, d\mu_g 
      \right)^\frac{p-1}{p} \!\!\!\|\xi\|_{L^p(\RR^n)} + \left( 
      \frac{1}{\epsilon_k^n} \int_M f'(u_k)u_k\, d\mu_g \right)^
      \frac{q-1}{q} \!\!\!\|\xi\|_{L^q(\RR^n)} \right]              \\
    &\leq C\left[\left(\frac{2\mu}{\mu-2}(m(J)+1)\right)^
      \frac{p-1}{p} \!\!\!\|\xi\|_{L^p(\RR^n)} + \left(\frac{2\mu}
      {\mu-2}(m(J)+1)\right)^\frac{q-1}{q} \!\!\!\|\xi\|_{L^q
      (\RR^n)} \right]
  \end{aligned}
  $$
  and this proves that $|E_{2,k}|$ tends to zero.
  Our second and last step is to prove that for any $\xi\in 
  C_0^\infty(\RR^n)$ $J'(w_k)(\xi)$ tends to $J'(w)(\xi)$ for $k$ tending to 
  infinity.  It is immediate that $\int_{\RR^n}\nabla w_k\cdot\nabla\xi\, 
  dz$ tends to $\int_{\RR^n}\nabla w\cdot\nabla\xi\, dz$.   By mean value 
  theorem there exists a function $\theta(z)$ with values in $(0,1)$ 
  such that
  $$
  \begin{aligned}
  \int_{\RR^n} &|f'(w_k(z))-f'(w(z))|\, |\xi(z)|\, dz \\
    &= \int_{\RR^n} |f''(\theta(z)w_k(z)+(1-\theta(z))
      w(z))|\, |w_k(z)-w(z)|\, |\xi(z)|\, dz\, .
  \end{aligned}
  $$
  By H\"older inequality the righthand side is bounded from above by
  $$
  \| w_k-w\|_{L^p(\Xi)} \|\xi\|_{L^p(\Xi)} \left( \int_{\RR^n} 
    |f''(\theta(z)w_k(z)+(1-\theta(z))w(z))|^\frac{p}{p-2} dz 
    \right)^\frac{p-2}{p}\, ,
  $$
  where $\| w_k-w\|_{L^p(\Xi)}$ tends to zero by $(i)$.  Besides we have
  $$
  \begin{aligned}
  \int_{\RR^n} &|f''(\theta(z)w_k(z)+(1-\theta(z))w(z))|^
      \frac{p}{p-2} dz                                       \\
    \leq &\; c_1\int_{\{z\in\Xi\;|\;|\theta(z)w_k(z)+(1-
      \theta(z))w(z)|\geq 1\}}  \hspace{-4cm}|\theta(z)
      w_k(z)+(1-\theta(z))w(z)|^p dz + c_1 \mbox{vol}\,(\Xi) \\
    \leq &\; c_1 2^{p-1} (\| w_k\|^p_{L^p(\Xi)}+
      \| w\|^p_{L^p(\Xi)}) + c_1 \mbox{vol}\,(\Xi)
  \end{aligned}
  $$
  and this quantity is bounded by a constant.

  \noindent $(iii)$  Let $t_k=t(w_k)$ be the multiplier defined in (ii), 
  Lemma \ref{lmm-tepsilonu}.  First of all we prove that there exist 
  $0<t_1\leq 1\leq t_2$ such that for all $k$ $t_1\leq t_k\leq t_2$.  
  Let $g_w(t)=J(tw)$.   By Lemma~\ref{lmm-tepsilonu} $(ii)$, it is 
  enough to find $t_1>0$ such that for all $t\in[0,t_1]$ 
  $g_{w_k}'(t)>0$ for all $k\in\NN$.  There holds
  $$
  \begin{aligned}
  g_{w_k}'(t) &= t \int_{\RR^n} |\nabla w_k(z)|^2 dz - \int_{\RR^n} 
      f'(tw_k(z))w_k(z)\, dz                                       \\
    &> t \int_{\RR^n} |\nabla w_k(z)|^2 dz - \frac{c_1t^{p-1}}{c_0} 
      \int_{\RR^n} f(w_k(z))\, dz\, .
  \end{aligned}
  $$
  Since we have
  $$
  \begin{aligned}
  \lefteqn{\int_{\RR^n} |\nabla w_k(z)|^2 dz \geq 
      \frac{h\epsilon^2_k}{H^\frac{n}{2}\epsilon^n_k} 
      \int_{B_g\left(x_k,\frac{R}{2}\right)} |\nabla u_k|^2_g
      \, d\mu_g}                                              \\
    &\geq \frac{2h}{H^\frac{n}{2}} \left( \frac{\epsilon^2_k}
      {2\epsilon^n_k} \int_{B_g\left(x_k,\frac{R}{2}\right)} 
      |\nabla u_k|^2_g\, d\mu_g - \frac{1}{\epsilon^n_k} 
      \int_{B_g\left(x_k,\frac{R}{2}\right)} f(u_k)\, d\mu_g 
      \right) \geq \frac{2h}{H^\frac{n}{2}} \gamma\, ,
  \end{aligned}
  $$
  where we have used the second equation of (\ref{gamma}), and
  $$
  \begin{aligned}
  \lefteqn{\int_{\RR^n} f(w_k(z))\, dz \leq \frac{1}{h^\frac{n}{2}
      \epsilon^n_k} \int_{B_g\left(x_k,\frac{R}{2}\right)} 
      f(u_k)\, d\mu_g}                                      \\
    &\leq \frac{2}{h^\frac{n}{2}(\mu-2)\epsilon^n_k} 
      \int_{B_g\left(x_k,\frac{R}{2}\right)} \left[ 
      \frac{1}{2} f'(u_k)u_k - f(u_k) \right] d\mu_g 
      \leq \frac{2(m(J)+1)}{h^\frac{n}{2}(\mu-2)} \, ,
  \end{aligned}
  $$
  then there exist $C_1,C_2>0$ such that $g_{w_k}'(t)>C_1t-C_2
  t^{p-1}$.  So we consider $t_1=\left(\frac{C_1}{C_2}
  \right)^\frac{1}{p-2}$.

  If $v$ is a function in the Nehari manifold $\mathcal{N}$, 
  $J(v)=\widetilde F_{\RR^n}(v)$, as defined 
  in (\ref{FtildeOmega}).   Then by property $(f1)$ $J(v)$ is 
  positive.    By Lemma~\ref{lmm-tepsilonu} $(ii)$, it is enough 
  to find $t_2>0$ such that for all $t\geq t_2$ $J(tw_k)<0$ 
  for all $k\in\NN$.  Since
  $$
  J(tw_k) = \frac{t^2}{2} \int_{\RR^n} |\nabla w_k(z)|^2 dz - 
    \int_{\RR^n} f(tw_k(z))\, dz
  $$
  and we already proved that $\{w_k\}_{k\in\NN}$ is bounded in 
  $\mathcal{D}^{1,2}(\RR^n)$, we still have to bound the second part 
  for $t\geq 1$
  $$
  \begin{aligned}
  \int_{\RR^n} f(tw_k(z))\, dz 
    &\geq c_0t^p \left( \int_{\{z\in\RR^n\;|\;|w_k(z)|\geq 1\}}
      \hspace{-1cm}|w_k(z)|^pdz + \int_{\{z\in\RR^n\;|\;|w_k(z)|
      \leq 1\}}\hspace{-1cm}|w_k(z)|^qdz \right)                 \\
    &> \frac{c_0t^p}{c_1} \int_{\RR^n} f''(w_k(z))(w_k(z))^2\, dz 
      > \frac{2c_0t^p}{c_1-2c_0} \widetilde F_{\RR^n}(w_k)       \\
    &\geq \frac{2c_0t^p}{(c_1-2c_0)H^\frac{n}{2}} \widetilde 
      F_{\epsilon_k,B_g\left(x_k,\frac{R}{2}\right)}(u_k) \geq 
      \frac{2c_0\gamma t^p}{(c_1-2c_0)H^\frac{n}{2}}\, ,
  \end{aligned}
  $$
  where we have used (\ref{gamma}).  So there exist $C_3,C_4>0$ such 
  that $J(tw_k)<C_3t^2-C_4t^p$ and $t_2=\left(\frac{C_3}{C_4}\right)^
  \frac{1}{p-2}$.

  By the boundedness of $t_k$ we conclude that up to subsequences 
  $t_k$ converges to $\bar t$ for $k$ tending to infinity.

  We apply the splitting lemma (Lemma \ref{lmm-splitting}) to the 
  sequence $t_kw_k$.  Then in the first case we have that 
  $t_kw_k$ converges strongly in $\mathcal{D}^{1,2}(\RR^n)$ to a 
  ground state solution $\bar w$.   It is easy to see that $t_kw_k$ 
  weakly converges to $\bar tw$, in fact for any $\xi\in C_0^\infty
  (\RR^n)$ there holds
  $$
  \begin{aligned}
  \left| \int_{\RR^n} \nabla(t_kw_k-\bar tw)\cdot\nabla\xi \right| 
    &= \left| \int_{\RR^n} \nabla(t_kw_k-\bar tw_k)\cdot\nabla\xi + 
      \int_{\RR^n} \nabla(\bar tw_k-\bar tw)\cdot\nabla\xi \right| \\
    &\leq |t_k-\bar t|\|\xi\|_{\mathcal{D}^{1,2}(\RR^n)}\|w_k\|_
      {\mathcal{D}^{1,2}(\RR^n)} + o(1) = o(1)\, .
  \end{aligned}
  $$
  We can conclude that $\bar w=\bar tw$.  In particular $w\not\equiv 0$ 
  and by the fact that both $\bar w$ and $w$ are in $\mathcal{N}$, 
  $\bar t=1$ and we have finished.

  Otherwise, there exist a sequence of points $\{y_k\}_{k\in\NN}$ 
  tending to infinity, 
  a ground state solution $U$ and a sequence of functions $\{w^0_k\}_
  {k\in\NN}$ such that, up to a subsequence $t_kw_k(z)=w^0_k(z)+
  U(z-y_k)$ for all $z\in\RR^n$ and $w^0_k$ tends strongly to zero.  
  We consider three different cases:  $\lim_{k\to\infty}|y_k|-
  \frac{R}{\epsilon_k}=2T>0$, $\lim_{k\to\infty}|y_k|-
  \frac{R}{\epsilon_k}=0$ and $\lim_{k\to\infty}\frac{R}{\epsilon_k}
  -|y_k|=2T>0$.   In the first case, since by definition $w_k\equiv 0$ 
  in $\RR^n\setminus B\left(0,\frac{R}{\epsilon_k}\right)$, 
  $w^0_k(z)=-U(z-y_k)$.  Then we have
  $$
  \begin{aligned}
  \lefteqn{\int_{\RR^n\setminus B\left(0,\frac{R}{\epsilon_k}
    \right)} |\nabla w^0_k(z)|^2 dz = \int_{\RR^n\setminus 
    B\left(0,\frac{R}{\epsilon_k}\right)} |\nabla U(z-y_k)|^2 dz} \\
    &\geq \int_{B(y_k,T)} |\nabla U(z-y_k)|^2 dz =  
      \int_{B(0,T)} |\nabla U(z)|^2 dz>0
  \end{aligned}
  $$
  and this is in contradiction with the fact that $w^0_k$ tends 
  strongly to zero.   If $\lim_{k\to\infty}|y_k|-
  \frac{R}{\epsilon_k}=0$, let $\pi(y_k)$ denote the 
  projection of $y_k$ onto the sphere centred in the origin with 
  radius $\frac{R}{\epsilon_k}$ and $T>0$.  Then
  $$
  \begin{aligned}
  \int_{\{z\in B(\pi(y_k),T)\;|\;|z|\geq\frac{R}
      {\epsilon_k}\}} 
    &|\nabla U(z-\pi(y_k))|^2 dz = \!\!\int_{\{z\in 
      B(0,T)\;|\;|z+\pi(y_k)|\geq\frac{R}{\epsilon_k}
      \}} \hspace{-0.4cm}|\nabla U(z)|^2 dz            \\
    &\geq \min_{\zeta\in S^n} \int_{\left\{z\in B(0,T)
      \;|\;z\cdot\zeta\geq 0\right\}} |\nabla U(z)|^2 
      dz = C>0
  \end{aligned}
  $$
  where $S^n$ is the unit sphere in $\RR^n$ and $z\cdot\zeta$ 
  is the scalar product in $\RR^n$.  Similarly to the first case 
  we have
  $$
  \begin{aligned}
  \int_{\RR^n\setminus B\left(0,\frac{R}{\epsilon_k}\right)} 
    |\nabla w^0_k(z)|^2 dz &= \int_{\RR^n\setminus B\left(0,
      \frac{R}{\epsilon_k}\right)} |\nabla U(z-y_k)|^2 dz    \\
    &\geq \int_{\left\{z\in B(y_k,T)\;|\;|z|\geq\frac{R}
      {\epsilon_k}\right\}} |\nabla U(z-y_k)|^2 dz           \\
    &= \int_{\left\{z\in B(\pi(y_k),T)\;|\;|z|\geq
      \frac{R}{\epsilon_k}\right\}} |\nabla U(z-
      \pi(y_k))|^2 dz + o(1)
  \end{aligned}
  $$
  and this is greater than $\frac{C}{2}$ for $k$ sufficiently large, 
  which is a contradiction.   Finally, if $\lim_{k\to\infty}
  \frac{R}{\epsilon_k}-|y_k|=2T>0$, for $k$ sufficiently large 
  $B(y_k,T)$ is contained in $B\left(0,\frac{R}{\epsilon_k}\right)$.   
  There holds
  $$
  \begin{aligned}
  \int_{B(y_k,T)}
    &\left[ \frac{1}{2}f'(U(z-y_k))U(z-y_k)-f(U(z-y_k)) \right] dz \\
    &= \int_{B(0,T)} \left[ \frac{1}{2}f'(U(z))U(z)-f(U(z)) \right] 
      dz = \gamma_0>0\, .
  \end{aligned}
  $$
  We consider the new sequence of points
  $$
  \tilde x_k = \exp_{x_k}(\epsilon_ky_k) \in B_g(x_k,R)\, .
  $$
  For any $k$ sufficiently large, let $U(\tilde x_k)$ be the 
  neighborhood of $\tilde x_k$ defined as $\exp_{x_k}(\epsilon_k
  B(y_k,T))$, then
  $$
  \begin{aligned}
  \lefteqn{\frac{1}{\epsilon^n_k} \int_{U(\tilde x_k)} \left[ 
      \frac{1}{2}f'(u_k)u_k-f(u_k) \right] d\mu_g}            \\ 
    &= \frac{1}{\epsilon^n_k} \int_{\epsilon_kB(y_k,T)} \left[ 
      \frac{1}{2}f'(u_k({\scriptstyle\exp_{x_k}(z)}))
      u_k({\scriptstyle\exp_{x_k}(z)})-f(u_k({\scriptstyle
      \exp_{x_k}(z)})) \right] |g_{x_k}(z)|^\frac{1}{2} dz    \\
    &\geq h^\frac{n}{2} \int_{B(y_k,T)} \left[ \frac{1}{2}f'
      (w_k(z))w_k(z)-f(w_k(z)) \right] dz\, .
  \end{aligned}
  $$
  Since $t_k\in(t_1,t_2)$ and using the properties of the function $f$ 
  we obtain
  $$
  \begin{aligned}
  \lefteqn{\int_{B(y_k,T)} \left[ \frac{1}{2}f'(w_k(z))w_k(z)-
      f(w_k(z)) \right] dz}                                    \\
    &\geq \int_{B(y_k,T)} \left[ \frac{1}{2}f' \left( 
      \frac{t_k}{t_2}w_k(z) \right) \frac{t_k}{t_2} w_k(z)-f 
      \left( \frac{t_k}{t_2}w_k(z) \right)\right] dz           \\
    &> \frac{(\mu-2)c_0}{(c_1-2c_0)t_2^q} \int_{B(y_k,T)} 
      \left[ \frac{1}{2}f'(t_kw_k(z)) t_kw_k(z)-f(t_kw_k(z)) 
      \right] dz\, .
  \end{aligned}
  $$
  By the splitting lemma we have
  $$
  \begin{aligned}
  \lefteqn{\int_{B(y_k,T)} \left[ \frac{1}{2}f'(t_kw_k(z)) 
      t_kw_k(z)-f(t_kw_k(z)) \right] dz}                    \\
    &= \int_{B(y_k,T)} \left[ \frac{1}{2}f'({\scriptstyle 
      w^0_k(z)+U(z-y_k)})({\scriptstyle w^0_k(z)+U(z-y_k)})
      -f({\scriptstyle w^0_k(z)+U(z-y_k)}) \right] dz       \\
    &= \int_{B(y_k,T)} \left[ \frac{1}{2}f'(U(z-y_k))(U(z-
      y_k))-f(U(z-y_k)) \right] dz + o(1)                   \\
    &=\gamma_0+o(1)\, .
  \end{aligned}
  $$
  So we have proved that for any $k$ sufficiently large
  \begin{equation}
  \frac{1}{\epsilon^n_k} \int_{U(\tilde x_k)} \left[ 
    \frac{1}{2}f'(u_k)u_k-f(u_k) \right] d\mu_g > \tilde\gamma_0 > 0\, .
  \label{tildegamma0}
  \end{equation}
  By definition, for $k$ big enough $U(\tilde x_k)$ is contained in 
  $B_g(\tilde x_k,R)$ and so we can substitute $x_k$ by $\tilde x_k$ 
  and $w_k$ by $\tilde w_k$, defined as in (\ref{w_k}) with the new 
  choice of points.  Steps $(i)$ and $(ii)$ are independent of $x_k$ 
  (provided $w_k$ is not identically zero) and so $\tilde w_k$ tends 
  weakly to a weak solution $\tilde w$.   It is possible to see that 
  there exists $\tilde T>0$ such that for any $k$ $U(\tilde x_k)
  \subset B_g(\tilde x_k,\epsilon_k\tilde T)$.  Then we have
  $$
  \begin{aligned}
  \int_{B(0,\tilde T)} &\left[ \frac{1}{2}f'(\tilde w_k(z))\tilde 
      w_k(z)-f(\tilde w_k(z)) \right] dz                         \\
    &\geq \frac{1}{H^\frac{n}{2}\epsilon_k^n} \int_{B_g(\tilde 
      x_k,\epsilon_k\tilde T)} \left[ \frac{1}{2}f'(u_k(x))u_k(x)
      -f(u_k(x)) \right] d\mu_g                                  \\
    &\geq \frac{1}{H^\frac{n}{2}\epsilon_k^n} \int_{U(\tilde 
      x_k)} \left[ \frac{1}{2}f'(u_k(x))u_k(x)-f(u_k(x)) \right] 
      d\mu_g\, .
  \end{aligned}
  $$
  By (\ref{tildegamma0}) and by the strong convergence of $\tilde 
  w_k$ to $\tilde w$ in $L^p(B(0,\tilde T))$, we conclude that
  $$
  \int_{B(0,\tilde T)} \left[ \frac{1}{2}f'(\tilde w(z))\tilde w(z)-
    f(\tilde w(z)) \right] dz \geq \frac{\tilde\gamma_0}{H^\frac{n}{2}}
  $$
  and so $\tilde w\not\equiv0$ and $\tilde w\in\mathcal{N}$.
  
  From now on we will write as before $w_k$ instead of $\tilde w_k$, 
  $x_k$ instead of $\tilde x_k$ and $w$ instead of $\tilde w$.  The 
  last step is to verify that $J(w)=m(J)$.  Let us consider the 
  following inequalities
  \begin{equation}
  \begin{aligned}
  m(J)+\delta_k &\geq J_{\epsilon_k}(u_k) = \frac{1}{\epsilon_k^n} 
      \int_M \left[ \frac{1}{2}f'(u_k)u_k - f(u_k) \right] d\mu_g  \\
    &\geq \int_{\RR^n} \left[ \frac{1}{2}f'(w_k)w_k - f(w_k) 
      \right] |g_{x_k}(\epsilon_k z)|^\frac{1}{2} dz\, .
  \end{aligned}
  \label{dismepsilon}
  \end{equation}
  We define the sequence of functions in $L^2(\RR^n)$:
  $$
  F_k(z) = \left[ \frac{1}{2}f'(w_k(z))w_k(z) - f(w_k(z)) 
    \right]^\frac{1}{2} |g_{x_k}(\epsilon_k z)|^\frac{1}{4}\, .
  $$
  By (\ref{dismepsilon}) this sequence is bounded in $L^2(\RR^n)$ and 
  there exists a weak limit $F\in L^2(\RR^n)$.  We prove that 
  \begin{equation}
  F(z) = \left[ \frac{1}{2}f'(w(z))w(z) - f(w(z)) \right]^\frac{1}{2}\, .
  \label{F(z)}
  \end{equation}
  Let $\xi$ be in $C_0^\infty(\RR^n)$.   On $\Xi$, the support of $\xi$, 
  $w_k$ strongly converges to $w$ in $L^p(\Xi)$.  So up to a subsequence 
  $w_k(z)$ converges to $w(z)$ almost everywhere.   Then pointwise
  $$
  F_k(z)\xi(z) \to \left[ \frac{1}{2}f'(w(z))w(z) - f(w(z)) 
    \right]^\frac{1}{2} \xi(z)
  $$
  for almost every $z\in\Xi$.  We can now apply Lebesgue theorem.  In 
  fact, there holds
  $$
  \begin{aligned}
  |F_k(z)|\,|\xi(z)| &< \left\{ \begin{array}{ll} 
      H^\frac{n}{4} \left( \frac{c_1}{2}-c_0 \right)^\frac{1}{2} 
      |w_k(z)|^\frac{p}{2} |\xi(z)| & \mbox{if }|w_k(z)|\geq 1\\
      H^\frac{n}{4} \left( \frac{c_1}{2}-c_0 \right)^\frac{1}{2} 
      |w_k(z)|^\frac{q}{2} |\xi(z)| & \mbox{if }|w_k(z)|\leq 1
      \end{array} \right.                                        \\
    &\leq H^\frac{n}{4} \left( \frac{c_1}{2}-c_0 \right)^
      \frac{1}{2} (1+|w_k(z)|^\frac{p}{2}) |\xi(z)|
  \end{aligned}
  $$
  and, since $w_k$ converges strongly to $w$ in $L^p(\Xi)$, there 
  exists 
  $W\in L^p(\Xi)$ such that for all $k$ $|w_k(z)|\leq W(z)$ almost 
  everywhere and $|F_k(z)|\, |\xi(z)|\leq H^\frac{n}{4} \left( 
  \frac{c_1}{2}-c_0 \right)^\frac{1}{2} (1+(W(z))^\frac{p}{2})
  |\xi(z)|\in L^2(\Xi)$.  So (\ref{F(z)}) is proved.
  By weak lower semicontinuity of the norm
  $$
  \| F\|^2_{L^2(\RR^n)} \leq \liminf_{k\to\infty} \| F_k\|^2_{L^2(\RR^n)}\, ,
  $$
  that is
  $$
  \int_{\RR^n} \left[ \frac{1}{2}f'(w)w - f(w) \right] dz \leq 
    \liminf_{k\to\infty} \int_{\RR^n} \left[ \frac{1}{2} f'(w_k)w_k 
    - f(w_k) \right] |g_{x_k}(\epsilon_k z)|^\frac{1}{2} dz\, .
  $$

  By this inequality and (\ref{dismepsilon}) we conclude 
  that
  $$
  \begin{aligned}
  m(J) &= \lim_{k\to\infty} m(J)+\delta_k \geq \lim_{k\to\infty} 
      J_{\epsilon_k}(u_k)                                       \\
    &\geq \liminf_{k\to\infty} \int_{\RR^n} \left[ \frac{1}{2}
      f'(w_k)w_k - f(w_k) \right] |g_{x_k}(\epsilon_k z)|^
      \frac{1}{2} dz                                            \\
    &\geq \int_{\RR^n} \left[ \frac{1}{2}f'(w)w-f(w) \right] dz 
      \geq m(J)\, .
  \end{aligned}
  $$
  
  $(iv)$ The equality is immediate from (\ref{dismepsilon}).
\end{proof}

We recall here Ekeland Principle (see for instance \cite{DeFigueiredo}).

\begin{definition}
  Let $X$ be a complete metric space and $\Psi:X\to\RR\cup
  \{+\infty\}$ be a lower semi-continuous function on $X$, bounded
  from below. Given $\eta>0$ and $\bar u\in X$ such that
  $$
  \Psi(\bar u) < \inf_{u\in X} \Psi(u) + \frac{\eta}{2}\, ,
  $$
  for all $\lambda >0$ there exists $u_\lambda\in X$ such that
  $$
  \Psi(u_\lambda) < \Psi(\bar u), \qquad d(u_\lambda,\bar u) <
    \lambda
  $$
  and for all $u\not= u_\lambda$ it holds
  $$
  \Psi(u_\lambda) < \Psi(u) + \frac{\eta}{\lambda} d(u_\lambda,u)\, .
  $$
\end{definition}

\begin{remark}
\label{rmk-importante}
  \begin{enumerate}
  \item We apply Lemma \ref{limSigma} when $u_k$ is a minimum 
    solution $u_k\in\mathcal{N}_{\epsilon_k}$, $J_{\epsilon_k}(u_k)
    =m_{\epsilon_k}$.  By $(iv)$ we have $\lim_{k\to\infty}m_{\epsilon_k}
    =m(J)$.  In particular for any $\delta>0$ there exists 
    $\epsilon_0=\epsilon_0(\delta)>0$ sufficiently small such that 
    for all $\epsilon\leq\epsilon_0$ $|m_\epsilon -m(J)|<\delta$.
  \item Applying Ekeland principle for $X=\Sigma_{\epsilon,\delta}$, 
    with $\epsilon\leq\epsilon_0(\delta)$ as in 1, we obtain that 
    for all $\bar u\in\Sigma_{\epsilon,\delta}$ there 
    exists $u_\delta\in\Sigma_{\epsilon,\delta}$ such that
    $$
    J_\epsilon(u_\delta) < J_\epsilon(\bar u), \qquad \frac{\epsilon}
      {\epsilon^\frac{n}{2}} \|u_\delta-\bar u\|_{H^1_2(M)} < 
      4\sqrt{\delta}
    $$
    and for all $u\in T\Sigma_{\epsilon,\delta}$
    \begin{equation}
    |J_\epsilon'(u_\delta)(u)| < \frac{\sqrt{\delta}\epsilon}
      {\epsilon^\frac{n}{2}} \| u\|_{H^1_2(M)}\, .
    \label{TSigma}
    \end{equation}
  \end{enumerate}
\end{remark}

\begin{proof}[Proof of Proposition~\ref{prp-x0}]
  We choose $\epsilon_0(\delta)$ as in point 1 of Remark 
  \ref{rmk-importante}.   We also assume that $\epsilon_0(\delta_0)$ 
  is less than $\epsilon_0$ in Lemma \ref{lmm-mepsilon>alpha}.
  
  By contradiction we assume that there is $\eta_0\in(0,1)$ such 
  that there exist two positive sequences $\{\delta_k\}_{k\in\NN}$, 
  $\{\epsilon_k\}_{k\in\NN}$ tending to zero as $k$ tends to 
  infinity and a sequence of functions $\{u_k\}_{k\in\NN}$, with 
  $u_k\in\Sigma_{\epsilon_k,\delta_k}$, and for any $x\in M$
  \begin{equation}
  \widetilde F_{\epsilon_k,B_g\left(x,\frac{r(M)}{2}\right)}(u_k) 
    \leq \eta_0\, m(J)\, .
  \label{proprieta}
  \end{equation}
  By Ekeland principle for any $k$ we can consider $\tilde u_k$ 
  as in 2 of Remark \ref{rmk-importante}.  Property 
  (\ref{proprieta}) becomes
  \begin{equation}
  \widetilde F_{\epsilon_k,B_g\left(x,\frac{r(M)}{2}\right)}
    (\tilde u_k) \leq \eta_1\, m(J)
  \label{proprietabis}
  \end{equation}
  with $\eta_1$ still in $(0,1)$.   To prove this we have to 
  evaluate the difference
  $$
  \frac{1}{\epsilon_k^n} \int_{B_g\left(x,\frac{r(M)}{2}\right)} 
    \left| \frac{1}{2}f'(\tilde u_k)\tilde u_k-f(\tilde u_k)-
    \frac{1}{2}f'(u_k)u_k+f(u_k) \right| d\mu_g\, ,
  $$
  which by mean value theorem can be written
  \begin{equation}
  \frac{1}{2\epsilon_k^n} \int_B \left| f''(u_k^*)u_k^*-f'
    (u_k^*) \right| |\tilde u_k-u_k|\, d\mu_g\, ,
  \label{star}
  \end{equation}
  where $B$ is $B_g\left(x,\frac{r(M)}{2}\right)$ and $u_k^*(x)
  =\theta(x)\tilde u_k(x)+(1-\theta(x))u_k(x)$ for a 
  suitable function $\theta(x)$ with values in $(0,1)$.    By 
  H\"older inequality (\ref{star}) is bounded from above by
  $$
  \frac{1}{2} \left( \frac{1}{\epsilon_k^n} \int_B |f''(u_k^*)
    u_k^*-f'(u_k^*)|^\frac{2n}{n+2}d\mu_g \right)^\frac{n+2}{2n} 
    \left( \frac{1}{\epsilon_k^n} \int_B |\tilde u_k-u_k|^
    \frac{2n}{n-2}d\mu_g \right)^\frac{n-2}{2n}.
  $$
  We prove that the first factor is bounded and the second one 
  is infinitesimal.  In fact, we have
  $$
  \begin{aligned}
  \left( \frac{1}{\epsilon_k^n} \int_B |\tilde u_k-u_k|^
      \frac{2n}{n-2} d\mu_g \right)^\frac{n-2}{2n} 
    &= \frac{\epsilon_k}{\epsilon_k^\frac{n}{2}} \|\tilde 
      u_k-u_k\|_{L^\frac{2n}{n-2}(B)}                     \\
    &\leq C\frac{\epsilon_k}{\epsilon_k^\frac{n}{2}} 
      \|\tilde u_k-u_k\|_{H^1_2(M)} < 4C\sqrt{\delta}\, .
  \end{aligned}
  $$
  The proof of the bound
  \begin{equation}
  \frac{1}{\epsilon_k^n} \int_B |f''(u_k^*)u_k^*-f'(u_k^*)|^
    \frac{2n}{n+2}d\mu_g \leq C
  \label{bound}
  \end{equation}
  for a positive constant $C$ can be found in the Appendix.

  We apply Lemma \ref{limSigma} to the sequences $\{\delta_k
  \}_{k\in\NN}$, $\{\epsilon_k\}_{k\in\NN}$ and $\{\tilde 
  u_k\}_{k\in\NN}$, obtaining a sequence of functions on $\RR^n$ 
  $\{w_k\}_{k\in\NN}$ (it is easy to see that (\ref{TSigma}) holds 
  for any $u\in H^1_2(M)$).   Let $w$ be the weak limit in $\mathcal{D}^
  {1,2}(\RR^n)$ of $w_k$.  Let $\eta_2$ be a constant in $(0,1)$ 
  such that $\eta_2>\frac{1+\eta_1}{2}$.   Since $J(w)=m(J)$, there 
  exists $T>0$ such that
  \begin{equation}
  \int_{B(0,T)} \left[ \frac{1}{2}f'(w(z))w(z)-f(w(z)) \right] dz 
    \geq \eta_2 m(J)\, .
  \label{dis2}
  \end{equation}
  On the other hand, up to a subsequence, we have
  \begin{equation}
  \begin{aligned}
  \lefteqn{\int_{B(0,T)} \left[ \frac{1}{2}f'(w)w-f(w) 
      \right] dz = \lim_{k\to\infty} \int_{B(0,T)} 
      \left[ \frac{1}{2}f'(w_k)w_k-f(w_k) \right] dz}   \\
    &= \lim_{k\to\infty} \frac{1}{\epsilon_k} \int_{B(0,
      \epsilon_k T)} \left[ \frac{1}{2}f'(\tilde u_k
      \circ\exp_{x_k})\tilde u_k\circ\exp_{x_k}-f(\tilde 
      u_k\circ\exp_{x_k}) \right] dz\, .
  \end{aligned}
  \label{dis2bis}
  \end{equation}
  By compactness the sequence $x_k$ converges (up to a subsequence) 
  to $\bar x$ and for any $z\in B(0,T)$ the limit of $|g_{x_k}
  (\epsilon_k z)|^\frac{1}{2}$ for $k$ tending to infinity is 
  $|g_{\bar x}(0)|^\frac{1}{2}=1$.   Since $\frac{2\eta_1}{1+\eta_1}
  \in(0,1)$, for $k$ sufficiently big for any $z\in B(0,\epsilon_k 
  T)$ we have $|g_{x_k}(z)|^\frac{1}{2}>\frac{2\eta_1}{1+\eta_1}$.  
  So the last limit in (\ref{dis2bis}) is less than
  $$
  \begin{aligned}
  \frac{1+\eta_1}{2\eta_1} &\lim_{k\to\infty} 
      \frac{1}{\epsilon_k} \!\int_{B(0,\epsilon_k T)} \!\!
      \left[ \frac{1}{2}f'(\tilde u_k\circ\exp_{x_k})\tilde 
      u_k\circ\exp_{x_k}-f(\tilde u_k\circ\exp_{x_k}) 
      \right]\! |g_{x_k}(z)|^\frac{1}{2}dz                 \\
    &=\frac{1+\eta_1}{2\eta_1} \lim_{k\to\infty} \frac{1}
      {\epsilon_k} \int_{B(x_k,\epsilon_k T)} \left[ 
      \frac{1}{2}f'(\tilde u_k)\tilde u_k-f(\tilde u_k) 
      \right] d\mu_g \leq \frac{1+\eta_1}{2} m(J)\, ,
  \end{aligned}
  $$
  where we have used property (\ref{proprietabis}).  By this 
  inequality together with (\ref{dis2bis}) and (\ref{dis2}) we get 
  $\eta_2\leq\frac{1+\eta_1}{2}$ wich is in contradiction with the 
  choice of $\eta_2$.
\end{proof}

It is now possible to prove the following proposition:

\begin{proposition}
\label{prp-beta}
  There exists $\delta_0\in(0,m(J))$ such that for any $\delta\in
  (0,\delta_0)$ there exists $\epsilon_0=\epsilon_0(\delta)>0$ and 
  for any $\epsilon\in(0,\epsilon_0)$ and $u\in\Sigma_{\epsilon,
  \delta}$ the barycentre $\beta(u)$ is in $M_{r(M)}$.
\end{proposition}

\begin{proof}
  By Proposition \ref{prp-x0}, for any $\eta\in(0,1)$ and for any 
  $u\in\Sigma_{\epsilon,\delta}$ with $\epsilon$ and $\delta$ 
  sufficiently small there exists a point $x_0$ such that
  $$
  \widetilde F_{\epsilon,B_g\left(x_0,\frac{r(M)}{2}\right)}(u) > 
    \eta\, m(J)\, .
  $$
  Since $u\in\Sigma_{\epsilon,\delta}$ we also have
  $$
  \widetilde F_{\epsilon,M}(u) \leq m(J)+\delta\, .
  $$
  We define
  $$
  \rho(u(x)) = \frac{\frac{1}{2}f'(u(x))u(x)-f(u(x))}{\int_M \left[ 
    \frac{1}{2}f'(u(x))u(x)-f(u(x)) \right] d\mu_g}\, .
  $$
  By the previous inequalities we have then
  $$
  \int_{B_g\left(x_0,\frac{r(M)}{2}\right)} \rho(u(x))\, d\mu_g > 
    \frac{\eta}{1+\frac{\delta}{m(J)}}\, .
  $$
  We can now esteem
  $$
  \begin{aligned}
  |\beta(u)-x_0| 
    &= \left| \int_M (x-x_0) \rho(u(x))\, d\mu_g \right|   \\
    &\leq \left| \int_{B_g\left(x_0,\frac{r(M)}{2}\right)} 
      \hspace{-1cm}(x-x_0) \rho(u(x))\, d\mu_g \right| + 
      \left| \int_{M\setminus B_g\left(x_0,\frac{r(M)}{2}
      \right)} \hspace{-1cm}(x-x_0) \rho(u(x))\, d\mu_g 
      \right|                                              \\
    &< \frac{r(M)}{2} + D \left( 1-\frac{\eta}{1+
      \frac{\delta}{m(J)}}\right)\, ,
  \end{aligned}
  $$
  where $D$ is the diameter of the manifold $M$.  For $\eta$ near 
  to $1$ and $\delta$ sufficiently small we obtain $\beta(u)\in 
  M_{r(M)}$.
\end{proof}

\section{The function $I_\epsilon$}

We prove now that the composition $I_\epsilon$ of $\phi_\epsilon$ 
and $\beta$ is well defined and homotopic to the identity on $M$:

\begin{proposition}
\label{prp-Iepsilon}
  There exists $\epsilon_0>0$ such that for any $\epsilon\in(0,
  \epsilon_0)$ the composition
  $$
  I_\epsilon=\beta\circ\phi_\epsilon:M\to M_{r(M)}
  $$
  is well defined and homotopic to the identity on $M$.
\end{proposition}

\begin{proof}
  Let us consider the function $H:[0,1]\times M\to M_{r(M)}$, defined 
  by $H(t,x)=tI_\epsilon(x)+(1-t)x$.  This function is a homotopy if 
  for any $t\in[0,1]$ $H(t,x)\in M_{r(M)}$.  It is enough to prove that 
  for any $x_0\in M$ $|I_\epsilon(x_0)-x_0|<r(M)$.   Since the support 
  of $\phi_\epsilon(x_0)$ is contained in $B_g(x_0,R)$
  $$
  \begin{aligned}
  I_\epsilon(x_0)-x_0 &= \int_M (x-x_0)\,\rho\left(\phi_\epsilon
      (x_0)(x)\right) d\mu_g = \int_{B_g(x_0,R)} \hspace{-1cm}
      (x-x_0)\,\rho\left(\phi_\epsilon(x_0)(x)\right) d\mu_g     \\
    &= \frac{\int_{B(0,R)} z\Phi(t_\epsilon(W_{x_0,\epsilon})
      W_{x_0,\epsilon}(\exp_{x_0}(z))) |g_{x_0}(z)|^\frac{1}{2} 
      dz}{\int_{B(0,R)} \Phi(t_\epsilon(W_{x_0,\epsilon})W_{x_0,
      \epsilon}(\exp_{x_0}(z))) |g_{x_0}(z)|^\frac{1}{2} dz}     \\
    &= \frac{\epsilon\int_{B\left(0,\frac{R}{\epsilon}\right)} 
      z\Phi(t_\epsilon(W_{x_0,\epsilon})W_{x_0,\epsilon}(\exp_
      {x_0}(\epsilon z))) |g_{x_0}(\epsilon z)|^\frac{1}{2} dz}
      {\int_{B\left(0,\frac{R}{\epsilon}\right)} \Phi(t_\epsilon
      (W_{x_0,\epsilon})W_{x_0,\epsilon}(\exp_{x_0}(\epsilon z))) 
      |g_{x_0}(\epsilon z)|^\frac{1}{2} dz}\, ,
  \end{aligned}
  $$
  where $\Phi$ is defined in (\ref{Phi}).   
  We recall that for any $\epsilon\in(0,1]$ and $x_0\in M$ 
  $t_1\leq t_\epsilon\left(W_{x_0,\epsilon}\right)\leq t_2$.  
  By definition of $\phi_\epsilon$, we have
  $$
  \int_{B\left(0,\frac{R}{\epsilon}\right)} \hspace{-0.8cm}\Phi
    (t_\epsilon(W_{x_0,\epsilon})W_{x_0,\epsilon}(\exp_{x_0}
    (\epsilon z))) |g_{x_0}(\epsilon z)|^\frac{1}{2} dz \geq 
    h^\frac{n}{2} \int_{B(0,R)} \hspace{-0.7cm}\Phi(t_1(U(z)-
    \widetilde U_R)) dz>0\, ,
  $$
  where $\widetilde U_R$ is the value $U(z)$ for any $z\in\RR^n$ 
  such that $|z|=R$.   Futhermore, we have
  $$
  \begin{aligned}
  \epsilon&\int_{B\left(0,\frac{R}{\epsilon}\right)} 
      \hspace{-0.8cm}|z| \Phi(t_\epsilon(W_{x_0,\epsilon})
      W_{x_0,\epsilon}(\exp_{x_0}(\epsilon z))) |g_{x_0}
      (\epsilon z)|^\frac{1}{2} dz \leq \epsilon H^
      \frac{n}{2} \int_{B\left(0,\frac{R}{\epsilon}\right)} 
      \hspace{-0.8cm}|z|\Phi(t_2U(z)) dz                   \\
    &< \frac{(c_1-2c_0)H^\frac{n}{2}\epsilon}{2} \left[ 
      \int_{\left\{z\in B\left(0,\frac{R}{\epsilon}\right)
      \;|\;t_2U(z)\geq 1\right\}} \hspace{-2cm}|z|t_2^p
      (U(z))^p dz + \int_{\left\{z\in B\left(0,\frac{R}
      {\epsilon}\right)\;|\;t_2U(z)\leq 1\right\}} 
      \hspace{-2cm}|z|t_2^q(U(z))^q dz \right].
  \end{aligned}
  $$
  Since $U$ is spherically symmetric and decreasing, there exists 
  $\rho_0>0$ such that the last quantity is equal to
  \begin{equation}
  \frac{(c_1-2c_0)H^\frac{n}{2}\epsilon}{2} \left[ 
    \int_{B(0,\rho_0)} \hspace{-0.5cm}|z|t_2^p(U(z))^p dz + 
    \int_{B\left(0,\frac{R}{\epsilon}\right)\setminus B(0,\rho_0)} 
    \hspace{-0.5cm}|z|t_2^q(U(z))^q dz \right].
  \label{quantita}
  \end{equation}
  Obviously, the integral
  $$
  \int_{B(0,\rho_0)} \hspace{-0.5cm}|z|t_2^p(U(z))^p dz \leq t_2^p
    \rho_0 \int_{B(0,\rho_0)} \hspace{-0.5cm}(U(z))^p dz
  $$
  is bounded.   For the second integral in (\ref{quantita}), we use 
  the well known inequality by Strauss (see \cite{Strauss}):
  $$
  \epsilon \int_{B\left(0,\frac{R}{\epsilon}\right)\setminus 
    B(0,\rho_0)} \hspace{-0.5cm}|z|(U(z))^q dz \leq C_n\|U\|_
    {\mathcal{D}^{1,2}(\RR^n)}^q \epsilon \int_{B\left(0,
    \frac{R}{\epsilon}\right)\setminus B(0,\rho_0)} \frac{|z|}
    {|z|^\frac{(n-2)q}{2}} dz
  $$
  where $C_n$ is a positive constant.  Then we conclude that 
  there exist two positive constants $C_1, C_2$ such that 
  (\ref{quantita}) is bounded from above by $C_1\epsilon+
  C_2\epsilon^\frac{(n-2)q-2n}{2}$, where the second exponent is 
  positive and so $|I_\epsilon(x_0)-x_0|$ tends to zero as 
  $\epsilon$ tends to zero.
\end{proof}

Finally, by standard arguments it is easy to see that the Palais-Smale 
condition holds for $J_\epsilon$ constrained on $\mathcal{N}_\epsilon$.

\section{The Morse theory result}

For an introduction to Morse theory we refer the reader to \cite{Milnor}, 
while for the applications to problems of functional analysis we mention 
\cite{B}.

Let $(X,Y)$ be a couple of topological spaces, with $Y\subset X$, and 
$H_k(X,Y)$ be the $k$-th homology group with coefficients in some 
field.  We recall the following definition:

\begin{definition}
  The Poincar\'e polynomial of $(X,Y)$ is the formal power series
  $$
  \mathcal{P}_t(X,Y) = \sum_{k=0}^\infty \dim[H_k(X,Y)] t^k\, .
  $$
  The Poincar\'e polynomial of $X$ is defined as $\mathcal{P}_t(X)=
  \mathcal{P}_t(X,\emptyset)$.
\end{definition}

If $X$ is a compact $n$-dimensional manifold $\dim[H_k(X)]$ is 
finite for any $k$ and $\dim[H_k(X)]=0$ for any $k>n$.  In particular 
$\mathcal{P}_t(X)$ is a polynomial and not a formal series.

We define now the Morse index.

\begin{definition}
  Let $J$ be a $C^2$ functional on a Banach space $X$ and let $u$ 
  be an isolated critical point of $J$ with $J(u)=c$.  The (polynomial) 
  Morse index of $u$ is defined as
  $$
  i_t(u) = \sum_{k=0}^\infty \dim[H_k(J^c,J^c\setminus\{u\})] t^k\, ,
  $$
  where $J^c=\{v\in X\;|\;J(v)\leq c\}$.  If $u$ is a non degenerate 
  critical point then $i_t(u)=t^{\mu(u)}$, where $\mu(u)$ is the 
  (numerical) Morse index of $u$ and represents the dimension of the 
  maximal subspace on which the bilinear form $J''(u)[\cdot,\cdot]$ 
  is negative definite.
\end{definition}

It is now possible to state Theorem \ref{trm-morse} more precisely:

\begin{theorem}
\label{trm-morse2}
  There exists $\epsilon_0>0$ such that for any $\epsilon\in(0,
  \epsilon_0)$, if the set $K_\epsilon$ of solutions of 
  equation (\ref{eq}) is discrete, then
  $$
  \sum_{u\in K_\epsilon} i_t(u) = t\mathcal{P}_t(M)+t^2[\mathcal{P}_t
    (M)-1]+t(1+t)\mathcal{Q}_\epsilon(t)\, ,
  $$
  where $\mathcal{Q}_\epsilon(t)$ is a polynomial with nonnegative 
  integer coefficients.
\end{theorem}

In the non-degenerate case, the above theorem becomes:

\begin{corollary}
\label{crl}
   There exists $\epsilon_0>0$ such that for any $\epsilon\in(0,
  \epsilon_0)$, if the set $K_\epsilon$ of solutions of 
  equation (\ref{eq}) is discrete and the solutions are non-degenerate, 
  then
  $$
  \sum_{u\in K_\epsilon} t^{\mu(u)} = t\mathcal{P}_t(M)+t^2[\mathcal{P}_t
    (M)-1]+t(1+t)\mathcal{Q}_\epsilon(t)\, ,
  $$
  where $\mathcal{Q}_\epsilon(t)$ is a polynomial with nonnegative 
  integer coefficients. 
\end{corollary}

Since we have proved that the composition $I_\epsilon$ of 
$\phi_\epsilon$ and $\beta$ from $M$ to $M_{r(M)}$ for $\epsilon$ 
sufficiently small is homotopic to the identity on $M$, the following 
equation holds (see \cite{BC}):
\begin{equation}
\mathcal{P}_t(\Sigma_{\epsilon,\delta})=\mathcal{P}_t(M) + 
  \mathcal{Z}(t)\, ,
\label{morse1}
\end{equation}
where $\mathcal{Z}(t)$ is a polynomial with nonnegative integer 
coefficients (here $\epsilon$ and $\delta$ are chosen as in 
Proposition \ref{prp-beta}).

Let $\alpha$ and $\epsilon$ be as in Lemma \ref{lmm-mepsilon>alpha}, 
$\delta>0$, then
\begin{equation}
\begin{aligned}
\mathcal{P}_t \left( J_\epsilon^{m(J)+\delta},J_\epsilon^
  \frac{\alpha}{2} \right) &= t \mathcal{P}_t(\Sigma_{\epsilon,
    \delta})\, ,                                                \\
\mathcal{P}_t \left( H^1_2(M),J_\epsilon^{m(J)+\delta} \right)
  &= t \left[ \mathcal{P}_t \left( J_\epsilon^{m(J)+\delta},
    J_\epsilon^\frac{\alpha}{2} \right) -t \right].
\end{aligned}
\label{morse2}
\end{equation}
By Morse theory we have
$$
\sum_{u\in K_\epsilon} i_t(u) = \mathcal{P}_t \left( H^1_2(M),
  J_\epsilon^{m(J)+\delta} \right) + \mathcal{P}_t \left( 
  J_\epsilon^{m(J)+\delta},J_\epsilon^\frac{\alpha}{2} \right) 
  + (1+t)\mathcal{Q}_\epsilon(t)\, ,
$$
where $\mathcal{Q}_\epsilon(t)$ is polynomial with nonnegative 
coefficients.   Using this relation with (\ref{morse1}) and 
(\ref{morse2}), we obtain Theorem \ref{trm-morse2} and Corollary 
\ref{crl}.  Theorem \ref{trm-morse} easily follows by evaluating 
the power series in $t=1$.

\section*{Appendix}

\begin{proof}[Proof of Lemma~\ref{lmm-estensione}]
  Given any $0<r<r(M)$, we can 
  choose $\rho<r$ small enough so that there exists a finite open 
  cover of $M_\rho$ $\{C_\alpha\}_{\alpha=1,
  \dots,k}$ of subsets of $\RR^N$ with smooth charts $\xi_\alpha:
  D_\alpha\subset\RR^N\to C_\alpha$ induced on $M_\rho$ by the 
  manifold structure of $M$.   We assume that $D_\alpha=Z_\alpha\times 
  T_\alpha$, with $Z_\alpha$ a subset of $\RR^n$ starshaped centred in 
  the origin and $T_\alpha$ the ball of $\RR^{N-n}$ with centre the 
  origin and radius $\rho$. For any $\alpha$ and any $(z,0)\in 
  Z_\alpha\times T_\alpha$, let $\xi_\alpha(z,0)\in\widetilde C_\alpha=
  C_\alpha\cap M$.  Viceversa for any $x\in\widetilde C_\alpha$, let 
  $\xi_\alpha^{-1}(x)=(z,0)$.

  We denote by $\{\psi_\alpha(y)\}_{\alpha=1,\dots,k}$ a partition of 
  unity subordinate to the cover $\{C_\alpha\}_{\alpha=1,\dots,k}$.   
  For all $y\in M_\rho$ we write $\xi_\alpha^{-1}(y)=
  (z_\alpha(y),t_\alpha(y))$.

  Given a function $u\in H^1_2(M)$, we define a function $v\in
  \mathcal{D}^{1,2}(M_r)$ by $v(y)\equiv 0$ for all $y\in M_r\setminus 
  M_\rho$ and
  $$
  v(y) = \sum_{\alpha=1}^k \psi_\alpha(y)\, u(\xi_\alpha(z_\alpha(y),0))
    \,\chi_\rho(|t_\alpha(y)|)
  $$
  for all $y\in M_\rho$, where $\chi_\rho$ is defined in (\ref{cut-off}).  

  \noindent\emph{Inequality (\ref{k1}).}  Let us write
  $$
  \begin{aligned}
  C_0 &= \left[ \sup_{i,j=1,\dots,N} \sup_{\alpha=1,\dots,k} 
        \sup_{y\in C_\alpha} \left( D_y(\xi_\alpha(z_\alpha(y),0)) 
        \right)_{ij} \right]^2\, ,                                \\
  C_1 &= \left[ \sup_{{i=1,\dots,N,\atop j=1,\dots, N-n}} 
        \sup_{\alpha=1,\dots,k} \sup_{y\in C_\alpha} \left( 
        D(t_\alpha(y)) \right)_{ij} \right]^2\, ,                 \\
  C_2 &= \sup_{\alpha=1,\dots,k} \sup_{y\in C_\alpha} (|\nabla
        \psi_\alpha(y)|^2+1)\, ,                                  \\
  C_3 &= \sup_{\alpha=1,\dots,k} \sup_{(z,t)\in D_\alpha} |\det 
        D(\xi_\alpha(z,t))|\, ,                                   \\
  C_4 &= \int_{\RR^{N-n}} \left[ (\chi_\rho(|t|))^2+(\chi_\rho'
        (|t|))^2 \right] dt\, .
  \end{aligned}
  $$
  Then we can estimate
  $$
  \begin{aligned}
  \int_{M_r} &|\nabla v(y)|^2 dy 
    \leq 2 \sum_{\alpha=1}^k \int_{C_\alpha} \left[ |\nabla
      \psi_\alpha(y)|^2 \left( u(\xi_\alpha(z_\alpha(y),0))\,
      \chi_\rho(|t_\alpha(y)|) \right)^2 \right.                    \\
    &\; + \left| \nabla_y \left( u(\xi_\alpha(z_\alpha(y),0)) 
      \right)\right|^2 (\psi_\alpha(y)\,\chi_\rho(|t_\alpha(y)|))^2 \\
    &\; + \left.\left| \nabla_y \left( \chi_\rho(|t_\alpha(y)|) 
      \right)\right|^2 (\psi_\alpha(y)\,u(\xi_\alpha(z_\alpha(y),
      0)))^2\right] dy                                              \\
    \leq &\; 2 \sum_{\alpha=1}^k \int_{C_\alpha} \left[ |\nabla
      \psi_\alpha(y)|^2 \left( u(\xi_\alpha(z_\alpha(y),0))\,
      \chi_\rho(|t_\alpha(y)|) \right)^2 \right.                    \\
    &\; + C_0 \left| \nabla u(\xi_\alpha(z_\alpha(y),0)) \right|^2 
      (\psi_\alpha(y)\,\chi_\rho(|t_\alpha(y)|))^2                  \\
    &\; + \left. C_1 \left( \chi_\rho'(|t_\alpha(y)|) \right)^2 
      (\psi_\alpha(y)\, u(\xi_\alpha(z_\alpha(y),0)))^2 \right] dy  \\
    \leq &\; \sum_{\alpha=1}^k \int_{C_\alpha} \left[ 2C_0 \left| 
      \nabla u(\xi_\alpha(z_\alpha(y),0)) \right|^2 (\chi_\rho
      (|t_\alpha(y)|))^2 \right.                                    \\
    &\; + \left. 2(1+C_1)C_2 (u(\xi_\alpha(z_\alpha(y),0)))^2 
      \left[ (\chi_\rho(|t_\alpha(y)|))^2+(\chi_\rho'(|t_\alpha
        (y)|))^2 \right] dy \right]                                 \\
    \leq &\; 2C_0C_3 \sum_{\alpha=1}^k \int_{D_\alpha} \left| \nabla 
      u(\xi_\alpha(z,0)) \right|^2 (\chi_\rho(|t|))^2\, dz\, dt     \\
    &\; + 2(1+C_1)C_2C_3 \sum_{\alpha=1}^k \int_{D_\alpha} (u
      (\xi_\alpha(z,0)))^2 \left[ (\chi_\rho(|t|))^2+(\chi_\rho'
      (|t|))^2 \right] dz\, dt                                      \\
    \leq &\; 2C_3(C_0+(1+C_1)C_2) \sum_{\alpha=1}^k \left[ 
      \int_{T_\alpha} (\chi_\rho(|t|))^2\, dt \int_{Z_\alpha} 
      \left| \nabla u(\xi_\alpha(z,0)) \right|^2\, dz \right.       \\
    &\; +\left. \int_{T_\alpha} \left[ (\chi_\rho(|t|))^2+
      (\chi_\rho'(|t|))^2 \right] dt \int_{Z_\alpha} (u(\xi_\alpha
      (z,0)))^2\, dz \right]                                        \\
    \leq &\; 2C_3(C_0+(1+C_1)C_2)C_4 \sum_{\alpha=1}^k 
      \int_{Z_\alpha} \left[ \left| \nabla u(\xi_\alpha(z,0)) 
      \right|^2 + (u(\xi_\alpha(z,0)))^2 \right] dz                 \\
    \leq &\; 2C_3(C_0+(1+C_1)C_2)C_4\frac{H}{h^\frac{n}{2}} 
      \sum_{\alpha=1}^k \int_{\widetilde C_\alpha} \left[ \left| 
      \nabla u(x) \right|^2_g + (u(x))^2 \right] d\mu_g\, .
  \end{aligned}
  $$
  One can easily see that there exists a constant $C_5>0$, depending only 
  on the charts $\xi_\alpha$ and on the partition of unity $\psi_\alpha$, 
  such that
  $$
  \sum_{\alpha=1}^k \int_{\widetilde C_\alpha} \left[ \left| 
    \nabla u(x) \right|^2_g + (u(x))^2 \right] d\mu_g \leq C_5 
    \| u\|^2_{H^1_2(M)}
  $$
  and by the Sobolev embedding of $H^1_2(M)$ in $L^2(M)$ (\ref{k1}) 
  is proved.

  \noindent\emph{Inequality (\ref{k2}).}  We show that for any 
  $s,t\in\RR$, $s+t\neq 0$
  $$
  f(s+t) > \frac{c_0\mu}{c_1}[f(s)+f(t)]\, .
  $$
  Let us consider first the case $|s+t|\geq 1$, $|s|\geq 1$ and 
  $|t|\geq 1$:
  $$
  f(s+t)\geq c_0 |s+t|^p\geq c_0(|s|^p+|t|^p)\geq\frac{c_0}{c_1}
    (f''(s)s^2+f''(t)t^2)>\frac{c_0\mu}{c_1}(f(s)+f(t))\, .
  $$
  If $|s+t|\geq 1$, $|s|\geq 1$ and $|t|<1$, we have:
  $$
  f(s+t)\geq c_0(|s|^p+|t|^p)\geq c_0(|s|^p+|t|^q)>\frac{c_0\mu}
    {c_1}(f(s)+f(t))\, .
  $$
  The same kind of inequalities holds true in the other cases.
  
  Hereafter, for all $y\in M_r$ we denote $v_\alpha(y)=\psi_\alpha(y) 
  u(\xi_\alpha(z_\alpha(y),0))\chi_\rho(|t_\alpha(y)|)$.  The 
  following integrals are always meant on the intersection with 
  the support of $v$:
  $$
  \begin{aligned}
  \int_{M_r} &f(v(y))\, dy = \int_{M_r} f\left( \sum_{\alpha=1}^k 
      v_\alpha(y) \right) dy > \frac{c_0\mu}{c_1} \sum_{\alpha=1}^k 
      \int_{C_\alpha} f(v_\alpha(y))\, dy                           \\
    &\geq \frac{c_0^2\mu}{c_1} \sum_{\alpha=1}^k \left[ 
      \int_{\{y\in C_\alpha\;|\;|v_\alpha(y)|\geq 1\}} \!\!
      |v_\alpha(y)|^p dy + \int_{\{y\in C_\alpha\;|\;|v_\alpha(y)|
      \leq 1\}} \!\!|v_\alpha(y)|^q dy \right]
  \end{aligned}
  $$
  For all $\alpha=1,\dots,k$ it is possible to choose $C_\alpha'
  \subset C_\alpha$ such that on this subset $\psi_\alpha(y)\geq
  \frac{1}{k}$.  Then the previous chain of inequalities is bounded from 
  below by
  \begin{equation}
  \begin{aligned}
  \frac{c_0^2\mu}{c_1k^q} &\sum_{\alpha=1}^k \left[ 
      \int_{\{y\in C_\alpha'\;|\;|v_\alpha(y)|\geq 1\}} |u(\xi_\alpha
      (z_\alpha(y),0))\,\chi_\rho(|t_\alpha(y)|)|^p dy \right.        \\
    &\; +\left. \int_{\{y\in C_\alpha'\;|\;|v_\alpha(y)|\leq 1\}} 
      |u(\xi_\alpha(z_\alpha(y),0))\,\chi_\rho(|t_\alpha(y)|)|^q dy 
      \right]
  \end{aligned}
  \label{dis}
  \end{equation}
  Let $D_\alpha'$ be the set $\xi_\alpha^{-1}(C_\alpha')$.  We consider 
  the following constants:
  $$
  \begin{aligned}
  C_5 &= \inf_{\alpha=1,\dots,k} \inf_{(z,t)\in D_\alpha} |\det 
        D(\xi_\alpha(z,t))|\, ,                                 \\
  C_6 &= \int_{\RR^{N-n}} (\chi_\rho(|t|))^q dt\, ,             \\
  C_7 &= \inf_{\alpha=1,\dots,k} \inf_{x\in\widetilde C_\alpha} 
        |\det D(z_\alpha(x))|\, .
  \end{aligned}
  $$
  The inequality (\ref{dis}) is bounded from below by
  $$
  \begin{aligned}
  \frac{c_0^2\mu C_5}{c_1k^q} &\sum_{\alpha=1}^k \left[ 
      \int_{\{(z,t)\in D_\alpha'\;|\;|v_\alpha(\xi_\alpha(z,t))|\geq 
      1\}} |u(\xi_\alpha(z,0))\,\chi_\rho(|t|)|^p dz\,dt \right.     \\
    &\; +\left. \int_{\{(z,t)\in D_\alpha'\;|\;|v_\alpha(\xi_\alpha
      (z,t))|\leq 1\}} |u(\xi_\alpha(z,0))\,\chi_\rho(|t|)|^q dz\, 
      dt \right]
  \end{aligned}
  $$
  $$
  \begin{aligned}
    \geq &\; \frac{c_0^2\mu C_5}{c_1k^q} \sum_{\alpha=1}^k \left[ 
      \int_{\{(z,t)\in D_\alpha'\;|\;|u(\xi_\alpha(z,0))|\geq 1\}} 
      |u(\xi_\alpha(z,0))|^p(\chi_\rho(|t|))^q dz\,dt \right.        \\
    &\; +\int_{\{(z,t)\in D_\alpha'\;|\;|u(\xi_\alpha(z,0))|\leq 1\}} 
      |u(\xi_\alpha(z,0))|^q(\chi_\rho(|t|))^q dz\, dt               \\
    &\; -\int_{\{(z,t)\in D_\alpha'\;|\;|v_\alpha(\xi_\alpha(z,t))|
      \leq 1,\,|u(\xi_\alpha(z,0))|\geq 1\}} |u(\xi_\alpha(z,0))|^p
      (\chi_\rho(|t|))^q dz\,dt                                      \\
    &\; +\left. \int_{\{(z,t)\in D_\alpha'\;|\;|v_\alpha(\xi_\alpha
      (z,t))|\leq 1,\,|u(\xi_\alpha(z,0))|\geq 1\}} |u(\xi_\alpha
      (z,0))|^q(\chi_\rho(|t|))^q dz\, dt \right]                    \\
    = &\; \frac{c_0^2\mu C_5C_6}{c_1k^q} \sum_{\alpha=1}^k \left[ 
      \int_{\{(z,0)\in D_\alpha'\;|\;|u(\xi_\alpha(z,0))|\geq 1\}} 
      |u(\xi_\alpha(z,0))|^p dz \right.                              \\
    &\; + \left. \int_{\{(z,0)\in D_\alpha'\;|\;|u(\xi_\alpha(z,0))|
      \leq 1\}} |u(\xi_\alpha(z,0))|^q dz \right]                    \\
    \geq &\; \frac{c_0^2\mu C_5C_6C_7}{c_1k^q} \sum_{\alpha=1}^k 
      \left[ \int_{\{x\in\widetilde C_\alpha\;|\; x\in C_\alpha',\, 
      |u(x)|\geq 1\}} |u(x)|^p dx \right.                            \\
    &\; + \left. \int_{\{x\in\widetilde C_\alpha\;|\; x\in C_\alpha',
      \, |u(x)|\leq 1\}} |u(x)|^q dx \right].
  \end{aligned}
  $$
  Since for all $x\in M$ the sum of the $\psi_\alpha(x)$ is one, there 
  exists $\hat\alpha$ such that $x\in C_\alpha'$.  Then for any $u\in 
  L^1(M)$
  $$
  \begin{aligned}
  \sum_{\alpha=1}^k \int_{C_\alpha'\cap M} |u(x)|\, dx 
    &= \sum_{\alpha=1}^k \int_{M} \chi_{C_\alpha'}(x)|u(x)|\, dx = 
      \int_{M} \left( \sum_{\alpha=1}^k\chi_{C_\alpha'}(x) \right) 
      |u(x)|\, dx                                                  \\
    &\geq \int_{M} |u(x)|\, dx\, .
  \end{aligned}
  $$
  This means that
  $$
  \begin{aligned}
  \sum_{\alpha=1}^k &\left[ \int_{\{x\in\widetilde C_\alpha\;|\; 
      x\in C_\alpha',\, |u(x)|\geq 1\}} |u(x)|^p dx + \int_{\{x
      \in\widetilde C_\alpha\;|\; x\in C_\alpha',\, |u(x)|\leq 1\}} 
      |u(x)|^q dx \right]                                           \\
    &\geq \int_{\{x\in M\;|\; |u(x)|\geq 1\}} |u(x)|^p dx + 
      \int_{\{x\in M\;|\; |u(x)|\leq 1\}} |u(x)|^q dx               \\
    &\geq \frac{1}{c_1} \int_M f''(u(x))(u(x))^2dx > \frac{\mu}{c_1} 
      \int_M f(u(x))\, dx \geq \frac{\mu}{c_1H^\frac{n}{2}} \int_M 
      f(u(x))\, d\mu_g\, .
  \end{aligned}
  $$

  \noindent\emph{Inequality (\ref{k3}).}  For $s>0$ 
  $f(s)$ is increasing.  Then we have
  $$
  \begin{aligned}
  \int_{M_r} f(v(y))\, dy &< \frac{c_1}{c_0\mu} \int_{M_r} 
      f(|v(y)|)\, dy \leq \frac{c_1}{c_0\mu} \int_{M_r} f \left( 
      \sum_{\alpha=1}^k |v_\alpha(y)| \right) dy                   \\
    &\leq \frac{c_1}{c_0\mu} \int_{M_r} f \left( \sum_{\alpha=1}
      ^k |\psi_\alpha(y)u(\xi_\alpha(z_\alpha(y),0))| \right) dy   \\
    &= \frac{c_1}{c_0\mu} \sum_{\beta=1}^k \int_{C_\beta} 
      \psi_\beta(y) f \left( \sum_{\alpha=1}^k |\psi_\alpha(y)
      u(\xi_\alpha(z_\alpha(y),0))| \right) dy                     \\
    &\leq \frac{c_1C_3}{c_0\mu} \sum_{\beta=1}^k \int_{D_\beta} 
      f \left( \sum_{\alpha=1}^k |\chi_{D_\alpha}(z,t)u(\xi_\alpha
      (z,0))| \right) dz\, dt                                      \\
    &\leq \frac{c_1C_3C_8}{c_0\mu} \sum_{\beta=1}^k 
      \int_{Z_\beta} f \left( \sum_{\alpha=1}^k |\chi_{Z_\alpha}(z)
      u(\xi_\alpha(z,0))| \right) dz\, ,
  \end{aligned}
  $$
  where $C_8$ is the volume of the ball of radius $\rho$ in $\RR^{N-n}$.  
  Proceeding with the chain of inequalities we obtain
  $$
  \begin{aligned}
  \sum_{\beta=1}^k &\int_{Z_\beta} f \left( \sum_{\alpha=1}^k 
      |\chi_{Z_\alpha}(z)u(\xi_\alpha(z,0))| \right) dz =  
      \sum_{\beta=1}^k \int_{\widetilde C_\beta} f \left( 
      \sum_{\alpha=1}^k |\chi_{\widetilde C_\alpha}(x)u(x)| 
      \right) dx                                             \\
    \leq &\;k \int_M f(k|u(x)|)\, dx                         \\
    < &\;\frac{kc_1}{\mu} \left[ \int_{\{x\in M\;|\;k|u(x)|
      \geq 1\}} k^p|u(x)|^pdx + \int_{\{x\in M\;|\;k|u(x)|
      \leq 1\}} k^q|u(x)|^qdx \right]                        \\
    = &\;\frac{kc_1}{\mu} \left[ \int_{\{x\in M\;|\;|u(x)|
      \geq 1\}} k^p|u(x)|^pdx + \int_{\{x\in M\;|\;|u(x)|
      \leq 1\}} k^q|u(x)|^qdx \right.                        \\
    &\; + \left. \int_{\{x\in M\;|\;|u(x)|\leq 1,\, k|u(x)|
      \geq 1\}} \hspace{-1cm}k^p|u(x)|^pdx - \int_{\{x\in M
      \;|\;|u(x)|\leq 1,\, k|u(x)|\geq 1\}} \hspace{-1cm}
      k^q|u(x)|^qdx \right]                                  \\
    \leq &\;\frac{kc_1}{\mu} \left[ \int_{\{x\in M\;|\;|u(x)|
      \geq 1\}} k^p|u(x)|^pdx + \int_{\{x\in M\;|\;|u(x)|
      \leq 1\}} k^q|u(x)|^qdx \right]                        \\
    \leq &\;\frac{k^{q+1}c_1}{c_0\mu} \int_M f(u(x))\, dx 
      \leq \frac{k^{q+1}c_1}{c_0\mu h^\frac{n}{2}} \int_M 
      f(u(x))\, d\mu_g\, .
  \end{aligned}
  $$

  \noindent\emph{Inequality (\ref{k4}).} The proof is analogous to 
  the proof of (\ref{k2}).
\end{proof}

We complete now the proof of Proposition~\ref{prp-x0}.

\begin{proof}[Proof of equation (\ref{bound})]
  The following inequalities hold:
  $$
  \begin{aligned}
  \frac{1}{\epsilon_k^n} &\int_B |f''(u_k^*)u_k^*-f'(u_k^*)|^
      \frac{2n}{n+2}d\mu_g                                    \\
    &\leq \frac{2^\frac{2n}{n+2}}{\epsilon_k^n} \int_B \left( 
      |f''(u_k^*)u_k^*|^\frac{2n}{n+2}+|f'(u_k^*)|^
      \frac{2n}{n+2} \right) d\mu_g                           \\
    &< \frac{2(2c_1)^\frac{2n}{n+2}}{\epsilon_k^n} \left( 
      \int_{\{x\in B\;|\;|u_k^*(x)|\geq 1\}} \hspace{-1cm}
      |u_k^*(x)|^\frac{(p-1)2n}{n+2} d\mu_g + \int_{\{x\in B
      \;|\;|u_k^*(x)|\leq 1\}} \hspace{-1cm}|u_k^*(x)|^
      \frac{(q-1)2n}{n+2} d\mu_g \right)                      \\
    &\leq \frac{2(2c_1)^\frac{2n}{n+2}}{\epsilon_k^n} \left( 
      \int_{\{x\in B\;|\;|u_k^*(x)|\geq 1\}} \hspace{-1cm}
      |u_k^*(x)|^p d\mu_g + \int_{\{x\in B\;|\;|u_k^*(x)|\leq 
      1\}} \hspace{-1cm}|u_k^*(x)|^q d\mu_g \right),
  \end{aligned}
  $$
  where in the last inequality we have used the fact that 
  $\frac{(p-1)2n}{n+2}<p$ and $\frac{(q-1)2n}{n+2}>q$.  We can 
  wright
  $$
  \begin{aligned}
  \int_{\{x\in B\;|\;|u_k^*(x)|\geq 1\}} \hspace{-1.8cm}
    &|u_k^*(x)|^p d\mu_g  + \int_{\{x\in B\;|\;|u_k^*(x)|\leq 
      1\}} \hspace{-1.8cm}|u_k^*(x)|^q d\mu_g                \\
    = &\int_{\{x\in B\;|\;|u_k^*(x)|\geq 1,\,|\tilde u_k(x)|
      \geq 1,\,|u_k(x)|\geq 1\}} \hspace{-1.8cm} |u_k^*(x)|^p 
      d\mu_g + \int_{\{x\in B\;|\;|u_k^*(x)|\geq 1,\,|\tilde 
      u_k(x)|\leq 1,\,|u_k(x)|\leq 1\}} \hspace{-1.8cm} 
      |u_k^*(x)|^p d\mu_g                                    \\
    &+ \int_{\{x\in B\;|\;|u_k^*(x)|\geq 1,\,|\tilde u_k(x)|
      \geq 1,\,|u_k(x)|\leq 1\}} \hspace{-1.9cm} |u_k^*(x)|^p 
      d\mu_g + \int_{\{x\in B\;|\;|u_k^*(x)|\geq 1,\,|\tilde 
      u_k(x)|\leq 1,\,|u_k(x)|\geq 1\}} \hspace{-1.9cm} 
      |u_k^*(x)|^p d\mu_g                                    \\
    &+ \int_{\{x\in B\;|\;|u_k^*(x)|\leq 1,\,|\tilde u_k(x)|
      \geq 1,\,|u_k(x)|\geq 1\}} \hspace{-1.9cm} |u_k^*(x)|^q 
      d\mu_g + \int_{\{x\in B\;|\;|u_k^*(x)|\leq 1,\,|\tilde 
      u_k(x)|\leq 1,\,|u_k(x)|\leq 1\}} \hspace{-1.9cm} 
      |u_k^*(x)|^q d\mu_g                                    \\
    &+ \int_{\{x\in B\;|\;|u_k^*(x)|\leq 1,\,|\tilde u_k(x)|
      \geq 1,\,|u_k(x)|\leq 1\}} \hspace{-1.9cm} |u_k^*(x)|^q 
      d\mu_g + \int_{\{x\in B\;|\;|u_k^*(x)|\leq 1,\,|\tilde 
      u_k(x)|\leq 1,\,|u_k(x)|\geq 1\}} \hspace{-1.9cm} 
      |u_k^*(x)|^q d\mu_g                                    \\
    \leq &\int_{\{x\in B\;|\;|u_k^*(x)|\geq 1,\,|\tilde u_k
      (x)|\geq 1,\,|u_k(x)|\geq 1\}} \hspace{-1.8cm} |u_k^*
      (x)|^p d\mu_g + \int_{\{x\in B\;|\;|u_k^*(x)|\geq 1,\,
      |\tilde u_k(x)|\leq 1,\,|u_k(x)|\leq 1\}} 
      \hspace{-1.8cm} |u_k^*(x)|^q d\mu_g                    \\
    &+ \int_{\{x\in B\;|\;|u_k^*(x)|\geq 1,\,|\tilde u_k(x)|
      \geq 1,\,|u_k(x)|\leq 1\}} \hspace{-1.9cm} |u_k^*(x)|^p 
      d\mu_g + \int_{\{x\in B\;|\;|u_k^*(x)|\geq 1,\,|\tilde 
      u_k(x)|\leq 1,\,|u_k(x)|\geq 1\}} \hspace{-1.9cm} 
      |u_k^*(x)|^p d\mu_g                                    \\
    &+ \int_{\{x\in B\;|\;|u_k^*(x)|\leq 1,\,|\tilde u_k(x)|
      \geq 1,\,|u_k(x)|\geq 1\}} \hspace{-1.9cm} |u_k^*(x)|^p 
      d\mu_g + \int_{\{x\in B\;|\;|u_k^*(x)|\leq 1,\,|\tilde 
      u_k(x)|\leq 1,\,|u_k(x)|\leq 1\}} \hspace{-1.9cm} 
      |u_k^*(x)|^q d\mu_g                                    \\
    &+ \int_{\{x\in B\;|\;|u_k^*(x)|\leq 1,\,|\tilde u_k(x)|
      \geq 1,\,|u_k(x)|\leq 1\}} \hspace{-1.9cm} |u_k^*(x)|^p 
      d\mu_g + \int_{\{x\in B\;|\;|u_k^*(x)|\leq 1,\,|\tilde 
      u_k(x)|\leq 1,\,|u_k(x)|\geq 1\}} \hspace{-1.9cm} 
      |u_k^*(x)|^p d\mu_g                                    \\
    \leq &\int_{\{x\in B\;|\;|u_k^*(x)|\geq 1,\,|\tilde u_k
      (x)|\geq 1,\,|u_k(x)|\geq 1\}} \hspace{-4.2cm} 2^p 
      \left( |\tilde u_k(x)|^p+|u_k(x)|^p \right) d\mu_g + 
      \int_{\{x\in B\;|\;|u_k^*(x)|\geq 1,\,|\tilde u_k
      (x)|\leq 1,\,|u_k(x)|\leq 1\}} \hspace{-4.2cm} 2^q 
      \left( |\tilde u_k(x)|^q+|u_k(x)|^q \right) d\mu_g     \\
    &+ \int_{\{x\in B\;|\;|u_k^*(x)|\geq 1,\,|\tilde u_k(x)|
      \geq 1,\,|u_k(x)|\leq 1\}} \hspace{-2.2cm} 2^p |\tilde 
      u_k(x)|^p d\mu_g + \!\int_{\{x\in B\;|\;|u_k^*(x)|\geq 1,
      \,|\tilde u_k(x)|\leq 1,\,|u_k(x)|\geq 1\}} 
      \hspace{-2.2cm} 2^p |u_k(x)|^p d\mu_g                  \\
    &+ \int_{\{x\in B\;|\;|u_k^*(x)|\leq 1,\,|\tilde u_k(x)|
      \geq 1,\,|u_k(x)|\geq 1\}} \hspace{-4.2cm} 2^p \left( 
      |\tilde u_k(x)|^p+|u_k(x)|^p \right) d\mu_g + \!
      \int_{\{x\in B\;|\;|u_k^*(x)|\leq 1,\,|\tilde u_k(x)|
      \leq 1,\,|u_k(x)|\leq 1\}} \hspace{-4.2cm} 2^q \left( 
      |\tilde u_k(x)|^q+|u_k(x)|^q \right) d\mu_g            \\
    &+ \int_{\{x\in B\;|\;|u_k^*(x)|\leq 1,\,|\tilde u_k(x)|
      \geq 1,\,|u_k(x)|\leq 1\}} \hspace{-2.2cm} 2^p |\tilde 
      u_k(x)|^p d\mu_g + \!\int_{\{x\in B\;|\;|u_k^*(x)|\leq 1,
      \,|\tilde u_k(x)|\leq 1,\,|u_k(x)|\geq 1\}} 
      \hspace{-2.2cm} 2^p |u_k(x)|^p d\mu_g                  \\
  \end{aligned}
  $$
  $$
  \begin{aligned}
  \leq &\int_{\{x\in B\;|\;|\tilde u_k(x)|\geq 1\}} 
      \hspace{-2.2cm} 2^p |\tilde u_k(x)|^p d\mu_g + 
      \int_{\{x\in B\;|\;|\tilde u_k(x)|\leq 1\}} 
      \hspace{-2.2cm} 2^q |\tilde u_k(x)|^q d\mu_g + \!
      \int_{\{x\in B\;|\;|u_k(x)|\geq 1\}} \hspace{-2.2cm} 
      2^p |u_k(x)|^p d\mu_g + \int_{\{x\in B\;|\;|u_k(x)|
      \leq 1\}} \hspace{-2.2cm} 2^q |u_k(x)|^q d\mu_g        \\
    \leq &\frac{2^q}{c_0} \int_M \left[ f(\tilde u_k) + 
      f(u_k) \right] d\mu_g\, .
  \end{aligned}
  $$
  Concluding there exists a constant $C>0$ such that
  $$
  \begin{aligned}
  \frac{1}{\epsilon_k^n} \int_B |f''(u_k^*)u_k^*-f'(u_k^*)|^
    \frac{2n}{n+2}d\mu_g &< \frac{C}{\epsilon_k^n} \int_M 
      \left[ f(\tilde u_k) + f(u_k) \right] d\mu_g           \\
    &\leq \frac{2C}{(\mu-2)} [J_{\epsilon_k}(\tilde u_k)+
      J_{\epsilon_k}(\tilde u_k)] \leq \frac{8Cm(J)}{(\mu-2)}
  \end{aligned}
  $$
  and this completes the proof of (\ref{bound}).
\end{proof}


\end{document}